\documentclass[11pt,a4paper]{article}
\usepackage{authblk}
\usepackage{indentfirst}
\usepackage{graphicx}
\usepackage{epsfig}

\usepackage{latexsym}
\usepackage{amsmath}
\usepackage{amssymb}
\usepackage{amsfonts}
\usepackage{mathrsfs}
\usepackage{pifont}

\usepackage{amsthm}
\usepackage{float}
\usepackage{subfigure}
\usepackage{color}
\usepackage{xcolor}
\usepackage{multicol}

\usepackage[colorlinks,linkcolor=blue,filecolor=black,citecolor=blue]{hyperref}

\newtheorem{theorem}{Theorem}[section]
\newtheorem{lemma}[theorem]{Lemma}
\newtheorem{prop}[theorem]{Proposition}
\newtheorem{coro}{Corollary}[theorem]

\newcommand{\mbf}{\boldsymbol}

\title{%Notes: Choice of learning samples for Helmholtz equation in domains with corners
An approximation for solutions to Helmholtz equation in polygons and its application in learning based numerical method
}

\author[1]{Yu Chen}
\author[2]{Jin Cheng \thanks{Corresponding author: jcheng@fudan.edu.cn}}
\affil[1]{{\small School of Mathematics, Shanghai University of Finance and Economics, Shanghai, China}}
\affil[2]{{\small School of Mathematical Sciences, Fudan University, Shanghai, China}}

\date{}
\begin{document}
\maketitle

\section*{Abstract}
A learning based method for Helmholtz equations was proposed in \cite{Chen2024}. Based on Runge's approximation, it is proved that when the solution can be continued outside the computational domain, using fundamental solutions as learning samples is sufficient to obtain the solution mapping with spectral accuracy. However, for complex geometries, such as domains with corners, approximation by fundamental solutions may only have logarithmic convergence due to the regularity issue, which is computationally impracticable. Moreover, approximation of a solution to the elliptic equation with those on larger domains is an ill-posed problem. In this work, we investigate how to feed additional learning samples to essentially get a practicable convergence rate and a stable algorithm, so as to ensure the generalizability of the learning based method. Error analysis and numerical examples illustrate flexibility and high-accuracy features.   
%\iffalse 

\section{Introduction}

%(A brief introduction for numerical method for Helmhotlz equation) 

The Helmholtz equation is a fundamental model for time-harmonic wave propagation in electromagnetics, acoustics, and seismology. Its accurate and efficient numerical computation is essential for wave simulation, optimal design, and inverse scattering, while also providing a means to investigate theoretical predictions. High-frequency computations remain particularly challenging: increasingly oscillatory solutions demand fine spatial resolution, numerical dispersion can accumulate into substantial phase errors, and the resulting indefinite linear systems are often difficult to solve efficiently. These difficulties are especially consequential for numerical investigations of increasing stability, whereby stability estimates for certain inverse problems improve as the frequency increases (see, e.g., \cite{Hrycak2004,Cheng2016}). Reliable numerical validation of such phenomenon therefore requires methods that control discretization errors and computational costs across a sufficiently broad range of frequencies.

A learning based method for Helmholtz equations is proposed in \cite{Chen2024}. The main idea is to reconstruct the solution mapping from various samples. One advantage of learning based numerical method is that solution to a new boundary input could be updated quickly once the solution operator is reconstructed. It is indicated that it is essential to employ PDE properties in the sample choice, reconstruction algorithm, and parameter choice, to ensure interpretability and generalizability. Based on quantitative Runge's approximation \cite{Ruland2018} and a priori choice of regularization parameter\cite{Cheng2000}, it was proved in \cite{Chen2024} that when the solution can be continued outside the computational domain, using fundamental solutions as learning samples is sufficient to obtain the solution mapping with spectral accuracy. However, if the solution cannot be continued to a larger domain, generally fundamental solutions may only have logarithmic convergence due to the regularity issue, which is computationally impracticable.

%(about sample choice) 

The importance of a suitable choice of learning samples becomes more prominent in complex geometries, such as in the presence of corners. For the Helmholtz equations, fundamental solutions with poles clustered at the corner are used to deal with the domain with corners \cite{Gopal2019}. It is revealed that the clustered pole solutions approximate the corner singularity. 

In two dimensional cases, the singularity caused by a corner for elliptic equations has been widely studied in, e.g.\cite{Grisvard1985}\cite{Kozlov1997}\cite{Kozlov2000}. Some special solutions to the equation are explicitly constructed\cite{Grisvard1985}. The weight of such solutions is usually determined approximately in FEM. These theoretical results suggest incorporating these special solutions into learning samples.

In analysis, one should further address issues concerning the denseness of sample space in solution space, the convergence, and the stability of solution operator reconstruction. In computations, denseness results along do not imply an acceptable convergence rate. What is the approximation property with the above samples, especially with complex geometry. Does a good approximation result imply the same performance for the learning method? Moreover, for flexibility, sample solutions in sufficiently larger domains are preferable to be adjusted to various domains. However, approximation with solutions on a larger domain is an ill-posed problem and raises the stability issues. These affect the numerical accuracy, efficiency, and generalizability, making a comprehensive analysis necessary.

In this work, a quantitative approximation for solutions to the Helmholtz equation in non-smooth domains is considered. Based on the quantitative Runge approximation and the singularity analysis with corners, the form of base functions and a computationally practical convergence rate are obtained.
The approximation provides a choice of learning samples for the reconstruction of solution operators. With suitable choice of regularization parameters, the error analysis is conducted, which essentially ensures interpretability and generalizability of the learning based method. Numerical examples illustrate the analysis.

%Propose an approximation for solutions in domains with corners with spectral convergence rate. Provide the error estimate for learning/approximation method. Regularization method is involved!

\section{Numerical method and main results}

Suppose that $u(x)$ satisfies Helmholtz equation with Dirichlet boundary condition
\begin{equation}%\label{problem1}
\begin{array}{ll}
	\Delta u(x)+k^2 u(x)=0, & \text{$x \in \Omega$}, \\
	u(x)=g(x), & \text{$x \in \partial\Omega $},\label{eq:Helm_boundary}
	\end{array}
\end{equation} 
where $\Omega$ is a bounded polygon in $\mathbb{R}^{2}$ and the boundary value is sufficiently smooth. $k>0$ is the wave number and $k^{2}$ is not a Dirichlet eigenvalue of $-\Delta$ in $\Omega$.

\paragraph{Basic algorithm}
Assuming there are sample solutions with the corresponding boundary values, we reconstruct the solution operator by these solution data.  Let $u_i(x),\ (i=1,2, \cdots, M)$ be a set of solutions to Helmholtz equation with different boundary values (training solutions), and let $T=\{x_j\}_{j=1}^N \subset\partial\Omega$ be the set of collocation points (interpolation nodes) on the boundary. Then generally for each sample solution $u_i(x)$, $x\in\Omega$, it can be approximated by 
$$
\sum_{j=1}^N A_j\left(x\right) u_i\left(x_j\right),\quad x_j\in T.
$$ 
Thus, there forms the following linear system approximately,
$$
\left(u_1(x), \cdots, u_M(x)\right)=\left(A_1(x), \cdots, A_N(x)\right)\left(\begin{array}{ccc}
	u_1\left(x_1\right) & \cdots & u_M\left(x_1\right) \\
	\vdots & \ddots & \vdots \\
	u_1\left(x_N\right) & \cdots & u_M\left(x_N\right)
\end{array}\right)\,,
$$
denoted by, 
$$
\mbf{b}_{1 \times M}=\mbf{a}_{1 \times N} V_{N \times M}.
$$
When $M>N$, it is an over-determined system. %To solve it we can take the pseudo inverse of $V$ and obtain
%$$
%\mbf{a}_{*}=\mbf{b} V^{\dagger}
%$$
%with proper singular value tolerance, or solving the Tikhonov minimization problem
We obtain the solution operator $a$ by solving the Tikhonov minimization problem,
\begin{equation}\label{eq:Tikhonov}
	\mbf{a}_{*}=\underset{\mbf{a} \in \mathbb{R}^{1\times N}}{\arg \min } \left\{\|\mbf{a}  V-\mbf{b}\|^{2}+\alpha \|\mbf{a}\|^2\right\},
\end{equation}
in which $\alpha\|\mbf{a}\|^2$ is the penalty term with an {\it a priori} given regularization parameter $\alpha>0$. Referring to the conclusions of Tikhonov regularization in \cite{Tikhonov1977}, \eqref{eq:Tikhonov} exists a unique minimum $\mbf{a}_{*}$ for all $\alpha>0$.

\paragraph{The equivalent approximation problem}

We will illustrate that in the above framework of the above learning method is equivalent to the approximation problem. The later in turn provide basis for sample choice and error estimate for LbMN.  

Let $Y$ be a Hilbert space that represents the solution (sub)space.%, e.g., $\mathbb{R}^M$, $\mathbf{\ell}^2$, 
\[
K_x: Y\rightarrow \mathbb{R},\quad
u(x)=K_x\mu \quad (=\sum_i\mu_i u_i(x)), \quad x\in \Omega.
\]

\[
K_f: Y\rightarrow L^2(\partial \Omega), \quad 
f(y)= K_f\, \mu\quad (=\sum_i\mu_i u_i(y)) ,\quad y\in\partial\Omega.
\]
Then, on one hand,
\[
\langle a,f\,\rangle=u(x)=\langle k_x, \mu\rangle_{Y},
\]
on the other hand,
\[
\langle a,f\,\rangle= \langle a, K_f\, \mu\rangle_{L^2(\partial\Omega)} =\langle K_f^* a,\mu\rangle_{Y}.
\]
The reconstruction problem is therefore to solve
\[
K_f^*a=k_x.
\]
Formally, it can be solved in the sense of pseudo inversion. As will be shown later, it is unstable and regularization should be introduced. Therefore, in real computations, we solve 
	\[
	a_*=\underset{\mbf{a} \in \mathbb{R}^{1\times N}}{\arg \min } \left\{\|K^{*}_f a - k_x\,\|_{Y}^{2}+\alpha\|a\|_{L^{2}(\partial\Omega)}^{2} \right\}.
	\]
The following proposition states that, in this framework, the learning method is equivalent to the approximation problem.
\begin{prop}
	If $a \in L^2(\partial\Omega)$ and $\mu \in \mathbb{R}^M$ are minimizers to
	\[
	\|K^{*}_f a - k_x\,\|_{Y}^{2}+\alpha\|a\|_{L^{2}(\partial\Omega)}^{2} \quad \text{and}\quad 
	\| K_f\,\mu - f\|^2_{L^2(\partial\Omega)}+\alpha \|\mu\|_Y^2
	\]
	respectively, then
	\[
	\langle a,f\,\rangle=\langle k_x,\mu\rangle \quad  (=\sum_{i} \mu_iu_i(x)).
	\]
\end{prop}	
	\;
	
	\begin{proof}
    By the linearity, $a$, and $\mu$ are solutions to
    \[
    (K_fK_f^*+\alpha I)a=K_f k_x; \quad (K^*_fK_f+\alpha I)\mu=K^*f
    \]
    respectively. Therefore,
        \[
        \langle a,f\rangle =\langle (K_fK_f^*+\alpha I)^{-1}K_fk_x,f \rangle=\langle k_x,K^*_f(K_fK_f^*+\alpha I)^{-*}f\rangle
        =\langle k_x,K^*_f(K_fK_f^*+\alpha I)^{-1}f\rangle,
        \]
        in which
        \begin{align*}
        K^*_f(K_fK_f^*+\alpha I)^{-1}& =(K^*_fK_f+\alpha I )^{-1}(K^*_fK_f+\alpha I )K^*_f(K_fK_f^*+\alpha I)^{-1}\\
        &=(K^*_fK_f+\alpha I )^{-1}K^*_f(K_fK_f^*+\alpha I)(K_fK_f^*+\alpha I)^{-1}=  (K^*_fK_f+\alpha I )^{-1}K^*_f.
        \end{align*}
        As a result,
        \[
        \langle a,f\rangle =\langle k_x, (K^*_fK_f+\alpha I )^{-1}K^*_f f \rangle =\langle k_x, \mu\rangle.
        \]
	\end{proof}

	%Recall:
	%\[
	%\langle\Phi_x,\mu\rangle = \langle a,f\,\rangle= \langle a, K_f\, \mu\rangle_{L^2(\partial\Omega)} =\langle K_f^* a,\mu\rangle_{Y}.
	%\]

%{\color{blue}
As a result, the error estimate of learning problem is equivalent to that of the approximation problem.%} 

\begin{figure}[H]
    \centering
    \includegraphics[width=0.3\linewidth]{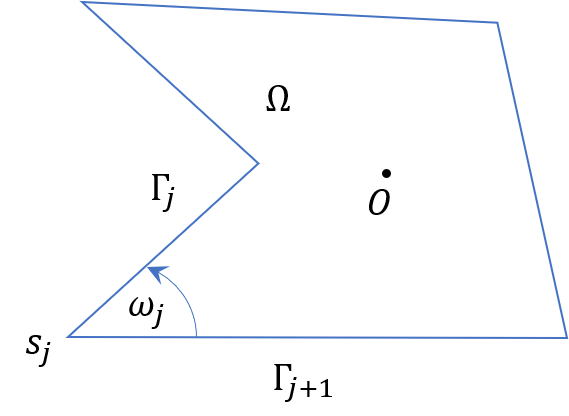}
    \caption{Notations of the polygon.}
    \label{polygon}
\end{figure}

\paragraph{Choice of learning samples and the main result}
We will prove that, for a star-shaped polygon, the basis 
\[
J_n(kr)e^{\mathrm{i}n\theta}, \quad n\in \mathbb{Z}
\]
centered at the origin (inside the polygon) and 
\begin{equation}\label{singularsolutions}
\sigma_{j,l}=
\begin{cases}
    J_{\lambda_{j,l}}(kr_j)\sin\left(\lambda_{j,l}\theta_j\right), & \lambda_{j,l}\notin\mathbb{Z}\\
   \frac{\partial J_{\lambda_{j,l}}}{\partial \lambda_{j,l}}\sin(\lambda_{j,l}\theta_j)+\theta_j J_{\lambda_{j,l}}(kr_j)\cos (\lambda_{j,l}\theta_j), & \lambda_{j,l}\in\mathbb{Z}
\end{cases}
,\quad j=1,2,\cdots, J; l=1,2,\cdots
\end{equation}
where $\lambda_{j,l}=l\pi/\omega_j$, centered at $j$th of $J$ corners form a complete set of learning samples for the Dirichlet boundary value problem and have a spectral convergence rate for smooth boundary values. Here $(r_j,\theta_j)$ is the local polar coordinates with origin at corner $s_j$ and $\theta_j=0$ on the open segment $\Gamma_{j+1}$ and $\theta_j=\omega_j$ on $\Gamma_j$ (see Fig. \ref{polygon}). 

For notations, set $I:=\{1,\cdots,J\}$, and denote
\[
X_{m,p_0}=\left\{g=(g_j)_{j=1}^J:
g_j\in W_{p_0}^{m+2-1/p_0}(\Gamma_j),\quad
g_j(s_j)=g_{j+1}(s_j)\right\},
\]
and
\[
\|g\|_{X_{m,p_0}}=
\left(\sum_{j=1}^J
\|g_j\|_{W_{p_0}^{m+2-1/p_0}(\Gamma_j)}^{p_0}\right)^{1/p_0}.
\]
%Indices at vertices are understood cyclically. This is a Banach space, since endpoint evaluation is continuous in the indicated edge spaces. It embeds continuously into $H^{1/2}(\partial\Omega)$.
Write
\[
\mathcal I_m=\{(j,l):j\in I, 0<\lambda_{j,l}\le m+1\},\qquad
|\mbf{c}|_1=\sum_{(j,l)\in\mathcal I_m}|c_{j,l}|.
\]

\paragraph{Main approximation result}
We will prove the following convergence results with the above choice of samples, which gives high precision with enough smooth boundary value.
%\begin{theorem}\label{main-approximation}
\begin{theorem}\label{main-approximation}
	%\color{teal}
	Let \(\Omega\Subset O_R\) be a strictly star-shaped polygon in $\mathbb{R}^2$. $u$ is a solution to the boundary value problem \eqref{eq:Helm_boundary} with $u|_{\partial\Omega}=g\in X_{m,p_0},\, p_0>2, m\in\mathbb N,\ m>1
	$.
	
	Then there exist coefficients $c_n$ and $c_{j,l}$, such that for $R>1$ with $k^2$ not being Dirichlet eigen value for $O_R$,
	\begin{equation}\label{solution-form}
	u_N(x)=\sum_{n=-N}^{N} \frac{c_n}{C_{R,n}}J_n(kr)e^{\mathrm{i}n\theta}  + \sum_{(j,l)\in \mathcal{I}_m} c_{j,l}\sigma_{j,l},\quad C_{R,n}=\frac{1}{|n|!}\left(\frac{kR}{2}\right)^{|n|}	
    \end{equation}
	converges to $u$ as $N\to \infty$, and with an error for $m>1$ as
	\[
	\|u-u_N\|_{L^2(\partial\Omega)} \leq  C \left(\frac{1}{\tilde{R}^{\sqrt{N}}}+\frac{1}{ N^{\frac{m}{2}}}\right)\|g\|_{X_{m,p_0}}.
	\]
	and
	\[
	\|u-u_N\|_{L^\infty(\partial\Omega)} \leq  C \left(e^{-\kappa_1\sqrt{N}}+\frac{1}{ N^{\kappa_2m}}\right)\|g\|_{X_{m,p_0}}.
	\]
	for some $\tilde{R}\in (1,R)$ and $\kappa_1,\kappa_2\in (0,1)$.
	%
	%
	%(\varepsilon+\frac{C}{R^N})  \|g\|_{H^{1/2}(\partial\Omega)}+ {\color{blue}C}\frac{\delta^{k-\frac{1}{2}-\eta}}{(k-2)!}\|g\|_{H^{k+3/2}(\partial \Omega)},\quad \forall \varepsilon, \eta>0.
	%\]
	Moreover, $c_n $ have bounds
	\[
	\sum_{n=-N}^N\sqrt{1+n^2} |c_n|^2 \leq C q^{-2N+c\sqrt{N}}\|g\|^2_{X_{m,p_0}},\quad |c_{j,l}|\leq C\|g\|_{X_{m,p_0}}, \quad    \lambda_{l,j} = \frac{l\pi}{\omega_j} \leq m+1. 
	\]
	for some $q\in (0,1)$ and $c>0$, and 
	\[
	q^{2N}\sum_{n=-N}^N\sqrt{1+n^2} |c_n|^2 \leq C e^{-a\sqrt{N}}\|g\|^2_{X_{m,p_0}}.
	\]

\end{theorem}

The estimate implies the following representation of a solution.

\begin{coro}
	With the assumptions in Theorem \ref{main-approximation}. If $g$ is piecewise smooth and satisfies the connection condition that $g_j(s_j)=g_{j+1}(s_j), j\in I$, then the solution can be expressed by
	\[
	u(x)=\sum_{n=-\infty}^{+\infty} \frac{C_n}{J_n(kR)}J_n(kr)e^{\mathrm{i}n\theta} + \sum_{\stackrel{j=1,\cdots,J}{l=1,\cdots }} c_{j,l}\sigma_{j,l}, \quad x\in \Omega. 
	\]
\end{coro}

With the approximation results and the estimate of the coefficients, we can further get the error analysis for the numerical method, which will be presented in Section 4.
\section{Analysis}
\subsection{Singularity decomposition}
%By the trace theorem for polygonal domains, one can transform the boundary value problem of the homogeneous equation into that of the nonhomogeneous equation with homogeneous boundary condition. However, for domains with corners, the source function may not lie in the image of spaces such as $\{u\in W_s^p(\Omega), u\}$ through the differential operator. Therefore, the shift theorem for smooth domains no longer holds and the codimension . 
\begin{lemma}[Grisvard 1985 \cite{Grisvard1985}]\label{reg-harmonic} %\color{teal}
Let $\Omega$ be a polygon, $f \in W_p^m(\Omega)$ and $g_j \in W_p^{m+2-1/p}(\Gamma_j), \quad j=1,\cdots,J$ be given. Assume that $( m\omega_j + 2\omega_j/q)/\pi$ is not an integer for any $j$ and that $g_j(S_j) = g_{j+1}(S_j)$. Then there exists $u$ which is a solution of the problem  
\[
\begin{cases}
\Delta u = f & \text{in } \Omega, \\
\gamma_j u = g_j & \text{on } \Gamma_j.
\end{cases}
\]
and constants $C_{j,l}$, such that
\[
u - \sum_{ 0 < \lambda_{j,l} < m+2/q} C_{j,l} \mathfrak{S}_{j,l} \in W_p^{m+2}(\Omega),
\]
where $\lambda_{j,l}=l\pi/\omega_j$ and
\[
\mathfrak{S}_{j,l}(r_j \mathrm{e}^{\mathrm{i}\theta_j}) = r_j^{\lambda_{j,l}} \sin\left(\lambda_{j,l} \theta_j \right) \eta_j(r_j \mathrm{e}^{\mathrm{i}\theta_j}),\quad \lambda_{j,l}\notin \mathbb{Z}
\]
when $\lambda_{j,l}$ is not an integer and
\[
\mathfrak{S}_{j,l}(r_j \mathrm{e}^{\mathrm{i}\theta_j}) = r_j^{\lambda_{j,l}} \bigl\{ \ln r_j \sin(\lambda_{j,l} \theta_j ) +\theta_j \cos(\lambda_{j,l} \theta_j ) \bigr\} \eta_j(r_j \mathrm{e}^{\mathrm{i}\theta_j}),\quad \lambda_{j,l}\in \mathbb{Z}
\]
when $\lambda_{j,l}$ is an integer. $\eta_j$ is some smooth cut-off function for the $j$th corner which equals 1 near the corner and $supp \,\eta_j \cap\Gamma_i=\emptyset, i\neq j,j+1$. 
\end{lemma}

\begin{theorem}\label{singularity-decomposition}
Assume that $u$ is the corresponding solution of $\Delta u + k^2u=0$ in the polygon $\Omega$ with boundary value $g\in W^{m+2-1/{p_0}}_{p_0}(\partial\Omega)$ for some $p_0>2$, $m\in \mathbb{N}_0$, and $k^2>0$ is not the Dirichlet eigen value of $-\Delta$. Assume that $g_j(s_j) = g_{j+1}(s_j)$ for $j\in I:=\{1,\cdots,J\}$. Then there exist numbers $c_{j,l}$ such that the solution to the boundary value problem has
%
%\vspace{-0.1in}
\[
u=u^r+\sum_{j\in I, \lambda_{j,l}\in (0,m+1]} c_{j,l}\sigma_{j,l} ,
\]
where $u^r\in H^{m+2}(\Omega)$, $\sigma_{j,l}$ is defined as \eqref{singularsolutions}. 
%
%\vspace{-0.2in}
%
%\[
%=J_{l\pi/\omega}(kr)\sin\left(\frac{l\pi\theta}{\omega}\right),\quad l=1,2,\cdots
%\] 
%$\omega$ is the inner angle of the corner and $l$ ranges over all integers such that $\pi/\omega\leq l\pi/\omega<2-2/p$.
\end{theorem}

\begin{proof}%\color{darkgray}
Choose $p\in (2,p_0)$ such that for any $d\in \{0,\cdots,m\}$, $( d\omega_j + 2\omega_j/q)/\pi$ is not an integer for any $j\in I:=\{1,\cdots,J\}$, where \(q=p/(p-1)\). 
    It holds that $g_j \in W_p^{m+2-1/p}(\Gamma_j), \quad j=1,\cdots,J$, and there exists $\hat{u}\in W^{m+2}_p(\Omega)$ such that $\gamma_j\hat{u}=g_j$ on $\Gamma_j$, by the trace theorem for a polygon of class $C^\infty$ (see, e.g., Theorem 1.5.2.8 in \cite{Grisvard1985}). According to standard well-posedness results for weak Dirichlet problems, there exists $\tilde{u}^{(1)}\in H^1(\Omega)$ such that
    \[
    \begin{cases}
    \Delta \tilde{u}^{(1)} + k^2 \tilde{u}^{(1)} = f & \text{in } \Omega, \\
    \gamma_j \tilde{u}^{(1)} = 0 & \text{on } \Gamma_j,\; j\in I.
    \end{cases}
    \]
    where $f=-(\Delta \hat{u} + k^2 \hat{u})\in L^p(\Omega)$ for $m=0$ and $W^m_p(\Omega)$ for $m\geq 1$. Then, $\tilde{u}^{(1)}=u-\hat{u}\in H^1_0(\Omega)$ satisfies
    \[
    \begin{cases}
    \Delta \tilde{u}^{(1)} = f-k^2\tilde{u}^{(1)}=:h_1 & \text{in } \Omega, \\
    \gamma_j \tilde{u}^{(1)} = 0 & \text{on } \Gamma_j,\; j\in I,
    \end{cases}
    \]
    where $h_1\in L^p(\Omega)$. Define $S_{j,l}=\sigma_{j,l}\eta_j$ where $\eta_j$ is the cut-off function for the $j$th corner as defined in Lemma \ref{reg-harmonic}. Then $\mathcal{L}S_{j,l}=0$ near the corner where $\mathcal{L}=\Delta+k^2I$. 
    By Lemma \ref{reg-harmonic}, there exist constants $\tilde{C}_{j,l}^{(1)}$ such that
    \[
    \tilde{u}^{(1)}-\sum_{ 0 < \lambda_{j,l} < 2/q} \tilde{C}_{j,l}^{(1)} \mathfrak{S}_{j,l} \in W_p^{2}(\Omega).
    \]
    %It can be shown by the asymptotic behavior of Bessel functions 
    After matching the leading coefficients, \(\mathfrak S_{j,l}-a_{j,l} S_{j,l}\in W_p^{d+2}(\Omega)\) for every term introduced at order \(d\). 
    Then there exist constants $C^{(1)}_{j,l}$ such that
    \[
    \sum_{0<\lambda_{j,l}<2/q}\tilde{C}^{(1)}_{j,l}\mathfrak{S}_{j,l}=\sum_{0<\lambda_{j,l}<2/q}C^{(1)}_{j,l}S_{j,l}+R,\quad R\in W_p^2(\Omega). 
    \]
    %Here $S_{j,l}=\eta_j\mathfrak{S}_{j,l}$ with $\eta_j$ is some smooth cut off function which is $1$ near the $j$th corner and with support has no intersection with other sides. 
    As a result,
    \[
    \tilde{u}^{(2)}:=\tilde{u}^{(1)}-\sum_{j\in I,\,0<\lambda_{j,l}<2/q}C^{(1)}_{j,l}S_{j,l}\in W_p^2(\Omega).
    \]
    When $m\geq 1$, since
    \begin{align*}
    \Delta\tilde{u}^{(2)}+k^2\tilde{u}^{(2)}&=(\Delta+k^2)\tilde{u}^{(1)}-\sum_{ 0 < \lambda_{j,l} < 2/q} C_{j,l}^{(1)} \mathcal{L}S_{j,l} \\
    &=f-\sum_{ 0 < \lambda_{j,l} < 2/q} C_{j,l}^{(1)} \mathcal{L}S_{j,l}=:f_2\in W_p^m(\Omega),
    \end{align*}
    and $S_{j,l}|_{\Gamma_i}\in C^\infty(\Gamma_i),S_{j,l}(s_j)=0,i=1,\cdots,J$ by its definition, $\tilde{u}^{(2)}$ satisfies
    \[
    \begin{cases}
    \Delta \tilde{u}^{(2)} = f_2-k^2\tilde{u}^{(2)}=:h_2 & \text{in } \Omega, \\
    \gamma_j \tilde{u}^{(2)} = g_2 & \text{on } \Gamma_j.
    \end{cases}    
    \]
    where $h_2\in W_p^1(\Omega)$, $g_{2,j}=g_2|_{\Gamma_j}\in W_p^{3-1/p}(\Gamma_j)$ and $g_{2,j}(s_j)=g_{2,j+1}(s_j)$ for $j\in I$. 
    
    Since \(\tilde u^{(2)}\in W_p^2(\Omega)\), the coefficients of the terms with \(\lambda_{j,l}<2/q\) in its new expansion vanish, by the linear independence of these corner terms modulo \(W_p^2(\Omega)\). 
    Then, applying Lemma \ref{reg-harmonic} yields that there exist constants $\tilde{C}^{(2)}_{j,l}$,
    \[
    \tilde{u}^{(2)}-\sum_{ 2/q < \lambda_{j,l} < 1+2/q} \tilde{C}^{(2)}_{j,l}\mathfrak{S}_{j,l}\in W_p^3(\Omega).\]
    Subsequently, since the normalized leading order of \(\mathfrak S_{j,l}-a_{j,l} S_{j,l}\) is \(r^{\lambda+2}\) for $\lambda>2/q$, there exist $C^{(2)}_{j,l}$ such that
    \[
    \tilde{u}^{(3)}:=\tilde{u}^{(2)}-\sum_{ 2/q < \lambda_{j,l} < 1+2/q} C^{(2)}_{j,l}S_{j,l}\in W_p^3(\Omega).
    \]
    Repeat the above processes until $W^{m+2}_p(\Omega)$, we have
    \[
    \tilde{u}^{(m+2)}:=\tilde{u}^{(m+1)}-\sum_{2/q+m-1<\lambda_{j,l}<2/q+m}C^{(m+1)}_{j,l}S_{j,l}\in W_p^{m+2}(\Omega).
    \]
    Finally, summing up the coefficients and considering that $(1-\eta_j)\sigma_{j,l}$ is smooth, we have that
    \begin{align*}
    u^r&=\hat{u}+\tilde{u}^{(1)}-\sum_{\lambda_{j,l}\in (0,m+2/q)}c_{j,l}\sigma_{j,l}\\
    &=\tilde{u}^{(1)}-\sum_{\lambda_{j,l}\in (0,m+2/q)}c_{j,l}S_{j,l} +\left(\hat{u}-\sum_{\lambda_{j,l}\in (0,m+2/q)}c_{j,l}(1-\eta_j)\sigma_{j,l}\right)   \in W_p^{m+2}\subset H^{m+2}.
    \end{align*}
    %where $S_{j,l}=(\eta-1)\sigma_{j,l}+\sigma_{j,l}$.
    %Considering the main singularity of $\mathfrak{S}_{j,l}$ and $\sigma_{j,l}$, we complete the proof.
    Absorb the terms $\sigma_{j,l}$ with respect to \(m+1<\lambda_{j,l}<m+2/q\) to $u^r$ since they belong to $H^{m+2}(\Omega)$, the proof is then complete.
\end{proof}

\begin{lemma}\label{regular-new}
With the assumptions of Theorem~\ref{singularity-decomposition}, let
$g\in X_{m,p_0}$, $p_0>2$. Then the decomposition
\[
u=u^r+\sum_{(j,l)\in\mathcal I_m}c_{j,l}\sigma_{j,l},
\qquad u^r\in H^{m+2}(\Omega),
\]
is unique and satisfies
\begin{equation}\label{alt-one-stability}
	\|u^r\|_{H^{m+2}(\Omega)}+
	\sum_{(j,l)\in\mathcal I_m}|c_{j,l}|
	\le C(m,p_0,k,\Omega)\|g\|_{X_{m,p_0}}.
\end{equation}
\end{lemma}

\begin{proof}
	Let $V_j=\Omega\cap B_\epsilon(s_j)$, with $\epsilon>0$ sufficiently
	small. By the asymptotic expansion of Bessel functions,
	\[
	\sigma_{j,l}=
	\begin{cases}
		b_\lambda r_j^\lambda\sin(\lambda\theta_j)
		+O(r_j^{\lambda+2}),&\lambda\notin\mathbb Z,\\[2pt]
		b_\lambda r_j^\lambda
		\{\log r_j\sin(\lambda\theta_j)+\theta_j\cos(\lambda\theta_j)\}
		+b_\lambda' r_j^\lambda\sin(\lambda\theta_j)
		+O(r_j^{\lambda+2}(1+|\log r_j|)),&\lambda\in\mathbb Z,
	\end{cases}
	\]
	where $\lambda=\lambda_{j,l}$ and
	$b_\lambda=(k/2)^\lambda/\Gamma(\lambda+1)\ne0$.
	For $\lambda\in\mathbb Z$, $r_j^\lambda\sin(\lambda\theta_j)$ is a
	polynomial. The singular leading terms satisfy
	\[
	\int_0^\epsilon r^{2(\lambda-m-2)}r\,\mathrm dr=\infty,
	\qquad 0<\lambda\le m+1,
	\]
	and hence $\sigma_{j,l}\notin H^{m+2}(V_j)$ for
	$(j,l)\in\mathcal I_m$.
	The smallest exponent with a nonzero coefficient cannot be cancelled
	by larger exponents; terms centered at other vertices are smooth in $V_j$.
	Thus
	\begin{equation}\label{alt-one-independence}
		\sum_{(j,l)\in\mathcal I_m}a_{j,l}\sigma_{j,l}
		\in H^{m+2}(\Omega)
		\quad\Longrightarrow\quad
		a_{j,l}=0,\qquad (j,l)\in\mathcal I_m.
	\end{equation}
	This proves uniqueness.
	
	Let $\mathcal P_k g=u$ be the Dirichlet solution operator. Since
	$k^2\notin\sigma(-\Delta_{\Omega,D})$,
	\[
	\|\mathcal P_k g\|_{H^1(\Omega)}
	\le C(k,\Omega)\|g\|_{H^{1/2}(\partial\Omega)}
	\le C\|g\|_{X_{m,p_0}}.
	\]
	By Theorem~\ref{singularity-decomposition} and
	\eqref{alt-one-independence}, the map
	\[
	T:X_{m,p_0}\longrightarrow Y_m:=
	H^{m+2}(\Omega)\times\mathbb C^{|\mathcal I_m|},
	\qquad Tg=(u^r,\mbf{c}),
	\]
	is well-defined and linear, where
	\[
	\|(v,\mbf{a})\|_{Y_m}
	:=\|v\|_{H^{m+2}(\Omega)}+
	\sum_{(j,l)\in\mathcal I_m}|a_{j,l}|.
	\]
	Assume that
	\[
	g_\nu\to g\quad\text{in }X_{m,p_0},\qquad
	Tg_\nu=(u_\nu^r,\mbf{c}_\nu)\to(v,\mbf{a})
	\quad\text{in }Y_m.
	\]
	Since $\sigma_{j,l}\in H^1(\Omega)$, it follows that
	\begin{align*}
		\mathcal P_k g_\nu&\rightarrow\mathcal P_k g
		&&\text{in }H^1(\Omega),\\
		\mathcal P_k g_\nu
		=u_\nu^r+\sum_{\mathcal I_m}c_{\nu,j,l}\sigma_{j,l}
		&\rightarrow v+\sum_{\mathcal I_m}a_{j,l}\sigma_{j,l}
		&&\text{in }H^1(\Omega).
	\end{align*}
	Therefore,
	\[
	\mathcal P_k g=v+\sum_{\mathcal I_m}a_{j,l}\sigma_{j,l},
	\qquad v\in H^{m+2}(\Omega)
	\quad\Longrightarrow\quad Tg=(v,\mbf{a}).
	\]
	The graph of $T$ is closed. By the closed graph theorem,
	\[
	\|u^r\|_{H^{m+2}(\Omega)}+|\mbf{c}|_1
	=\|Tg\|_{Y_m}\le C(m,p_0,k,\Omega)\|g\|_{X_{m,p_0}},
	\]
	which completes the proof.
\end{proof}

\subsection{Approximation of the regular part}

\begin{lemma}[Burenkov and Gorbunov 1997]\label{norm-extension}
	For any bounded set $\Omega\subset \mathbb{R}^2$ with boundary of class Lip 1, there exist $h(\Omega)>0$ and an extension operator $E: W^{s}_p(\Omega)\rightarrow W^{s}_p(\Omega_h)$,  $u\mapsto \tilde{u}=Eu$, such that $E$ is continuous with 
	\[
	\|E\|_{W^{\,s}_p(\Omega)\rightarrow W_p^{\,s} (\Omega_h)}\leq C(\Omega)^{\,s}, \quad s\in \mathbb{N}, 1\leq p\leq +\infty.
	\]
\end{lemma}

Rm: this result to make the constant rely explicitly on the order $s$

\begin{lemma}\label{source-extension}
	Let $\Omega\subset\mathbb R^2$ be a bounded polygonal domain, and set
	$\Omega_\delta=\{x:\operatorname{dist}(x,\Omega)<\delta\}$ and
	$A_\delta=\Omega_\delta\setminus\overline\Omega$. %	
	Let $\Omega_h$ be a fixed domain with $\Omega\Subset\Omega_h$, and let $m\in\mathbb N_0$. If $v\in H^{m}(\Omega_h)$ and $v=0$ a.e. on $\Omega$, then with
	$\omega(\delta)=\|v\|_{L^2(A_\delta)}$, one has
	\begin{equation}\label{extension-dependency}
		\lim_{\delta\to 0}\omega(\delta)=0\quad(m=0),
		\qquad
		\omega(\delta)\leq (C_\Omega\delta)^m
		\|v\|_{H^m(A_\delta)}\quad(m\geq1),
	\end{equation}
	for every $\delta>0$ such that $\Omega_\delta\Subset\Omega_h$, where $C_\Omega$ is independent of $m$, $v$, and $\delta$; the limit is taken with $\Omega_h$ and $v$ fixed.
\end{lemma}

\begin{proof}
	Choose $r_\Omega>0$ such that $\Omega\cap B(s_i,r_\Omega)$ is a sector for every vertex $s_i$, and the balls $B(s_i,r_\Omega)\cap B(s_j,r_\Omega)=\emptyset, i\neq j$. Let $0<\delta\leq r_\Omega$.
	For a.e. $x\in A_\delta$, let $\pi(x)$ be its unique nearest point in $\overline\Omega$:
	$|x-\pi(x)|=\operatorname{dist}(x,\overline\Omega)$.
	For an edge $\Gamma_i$ which is an open segment, set
	\[
	E_i=\{x\in A_\delta:\pi(x)\in \Gamma_i\},\qquad
	V_i=\{x\in A_\delta:\pi(x)=s_i\}.
	\]
	These sets partition $A_\delta$ up to a null set.
	Then, $E_i$ can be represented as %For an edge $e$ with outward normal $n_e$,
	\[
	E_i= \{y+t \nu_i :y\in \Gamma_i,\ 0<t<L_e(y)\},\qquad
	0\leq L_e(y)\leq\delta.
	\]
	where $\nu_i$ is the outward normal on $\Gamma_i$. %	The normal fibers are intervals. 
	For a.e. $y\in \Gamma_i$, by slicing property of Sobolev
	functions and the zero traces of $D^\alpha v$ ($|\alpha|<m$),
	\[
	g_y(t):=v(y+t \nu_i )\in H^m(0,L_e(y)),\qquad
	g_y^{(j)}(0)=0\quad(0\leq j<m).
	\]
	The integral Taylor formula, valid a.e. for one-dimensional
	$H^m$ functions, and Young's inequality yield
	\begin{align}
		g_y(t)
		&=\frac{1}{(m-1)!}\int_0^t(t-s)^{m-1}g_y^{(m)}(s)\,ds,
		\notag\\
		\|g_y\|_{L^2(0,L_e(y))}
		&\leq\frac{L_e(y)^m}{m!}
		\|g_y^{(m)}\|_{L^2(0,L_e(y))}
		\leq\frac{\delta^m}{m!}
		\|\partial_{\nu_i}^m v(y+\cdot \nu_i)\|_{L^2(0,L_e(y))}.
		\label{normal-taylor-bound}
	\end{align}
	%The normal-coordinate Jacobian is $1$; 
	Therefore, 
	\begin{equation}\label{edge-layer-bound}
		\|v\|_{L^2(E_i)}
		\leq \frac{\delta^m}{m!}\|\partial_{\nu_i}^m v\|_{L^2(E_i)}
		\leq\frac{(\sqrt2\,\delta)^m}{m!}\|v\|_{H^m(E_i)}.
	\end{equation}
	
	For a vertex $s_i$, $V_i\subset B(s_i,\delta)$. Write
	$e_\theta=(\cos\theta,\sin\theta)$, and let $I_i$ be the
	exterior angular interval in $B(s_i,r_\Omega)$, with length
	$\beta_i\leq2\pi$. Since $D^\alpha v=0$ in the interior sector,
	the angular traces of $D^\alpha v(s_i+r e_\theta)$ vanish at
	$\partial I_i$ for $|\alpha|<m$ and a.e. $0<r<\delta$.
	For $h\in H_0^1(I_i)$, the one-dimensional Dirichlet
	Poincar\'e inequality gives
	\[
	\int_{I_i}|h|^2\,d\theta
	\leq\frac{\beta_i^2}{\pi^2}
	\int_{I_i}|h'|^2\,\mathrm{d}\theta.
	\]
	Set $\alpha!=\alpha_1!\alpha_2!$ and
	\[
	F_k(r)=\sum_{|\alpha|=k}\frac{k!}{\alpha!}
	\int_{I_i}|D^\alpha v(s_i+r e_\theta)|^2\,\mathrm{d}\theta.
	\]
	Since
	\[
	\sum_{|\alpha|=k}\frac{k!}{\alpha!}|\nabla D^\alpha v|^2
	=\sum_{|\gamma|=k+1}\frac{(k+1)!}{\gamma!}|D^\gamma v|^2,
	\qquad |\partial_\theta D^\alpha v|\leq r|\nabla D^\alpha v|,
	\]
	\[
	F_k(r)\leq\frac{\beta_i^2r^2}{\pi^2}F_{k+1}(r)
	\leq4r^2F_{k+1}(r),\qquad 0\leq k<m.
	\]
	Since $\max_{0\leq j\leq m}\binom mj\leq2^{m-1}$ for $m\geq1$,
	\[
	F_m(r)\leq2^{m-1}\sum_{|\alpha|=m}
	\int_{I_i}|D^\alpha v(s_i+r e_\theta)|^2\,d\theta.
	\]
	Consequently,
	\begin{align}
		\|v\|_{L^2(V_i)}^2
		&\leq\int_0^\delta F_0(r)r\,dr
		\leq(4\delta^2)^m\int_0^\delta F_m(r)r\,dr
		\notag\\
		&\leq2^{m-1}(4\delta^2)^m
		\|v\|_{H^m(B(s_i,\delta)\cap A_\delta)}^2.
		\label{vertex-layer-bound}
	\end{align}
	Summing over disjoint $E_i$ and pairwise disjoint vertex balls,
	\[
	\|v\|_{L^2(A_\delta)}^2
	\leq\left(\frac{2^m}{(m!)^2}
	+\tfrac12\,8^m\right)\delta^{2m}
	\|v\|_{H^m(A_\delta)}^2
	\leq(2\sqrt2\,\delta)^{2m}
	\|v\|_{H^m(A_\delta)}^2.
	\]
	Here $2^m/(m!)^2\leq\tfrac12\,8^m$ for $m\geq1$.
	For $\delta>r_\Omega$ with $\Omega_\delta\Subset\Omega_h$,
	\[
	\|v\|_{L^2(A_\delta)}
	\leq\|v\|_{H^m(A_\delta)}
	\leq(\delta/r_\Omega)^m\|v\|_{H^m(A_\delta)}.
	\]
	Thus $C_\Omega=\max\{2\sqrt2,r_\Omega^{-1}\}$.
	For $m=0$, $|A_\delta|\downarrow 0$ and $v\in L^2(\Omega_h)$ imply $\omega(\delta)\to 0$ by absolute integrability.
\end{proof}

\begin{lemma}[Approximation of the regular part $u_r$]\label{harmonic-approx2}
	Let $\Omega\subset \mathbb{R}^2$ be a bounded polygon domain. If $u_r\in H^{m+2}(\Omega)$ satisfies Helmholtz equation in $\Omega$, % and  $k^2$ is not a Dirichlet eigen value of $-\Delta$ in $\Omega$, $m\in\mathbb{N}$. 
    Then for $\forall \epsilon>0$ and $m=0$, there exist a bounded smooth domain $\hat{\Omega}$ such that $\Omega\Subset\hat{\Omega}$ and a solution $\hat{u}_r$ of the Helmholtz equation in $\hat{\Omega}$ such that 
	\[
	\|\hat{u}_r|_{\Omega}-u_r\|_{H^1(\Omega)} \leq \epsilon.
	\]
	
	Moreover, for $m\geq 1$, there exists $h(\Omega,k)>0$ such that for $0<\delta<h(\Omega)$, one can choose a bounded smooth domain $\hat{\Omega}$ with $\Omega_{\delta/2}\Subset \hat{\Omega}\Subset \Omega_{\delta}$ and
    % and let $\hat{\Omega}\Supset \Omega$ be a domain with smooth boundary satisfying , then
    $\hat{u}^\delta_r\in H^1(\hat{\Omega})$ satisfying Helmholtz equation in $\hat{\Omega}$ such that
	\begin{equation}
		\|\hat{u}^\delta_r|_{\Omega}-u_r\|_{H^1(\Omega)} \leq 
		  C(m,\Omega,k)\delta^m \|u_r\|_{H^{m+2}(\Omega)},\quad \|\hat{u}_r^\delta\|_{H^{m+2}(\hat{\Omega})}
	\leq C(\Omega,m,k)\|u_r\|_{H^{m+2}(\Omega)}.
	\end{equation}
\end{lemma}

\begin{proof}

        By Lemma \ref{norm-extension} with $s=m+2$, $p=2$, there exists $h(\Omega)>0$ and an extension of $u_r$: $\tilde{u}_r=Eu_r$, such that        
       % since $u_r\in H^{m+2}(\Omega)$, and an continuous extension operator $E: H^{m+2}(\Omega)\rightarrow H^{m+2}(\Omega_h)$,  , s.t.
		\[
		\|Eu_r\|_{H^{m+2}(\Omega_h)}\leq C(\Omega)^{m+2} \|u_r\|_{H^{m+2}(\Omega)}.
		\]
        where $\Omega_h=\{x:dist(x,\Omega)<h\}$.

		%\item Extension of source============ 
		
		Define $f\triangleq \mathcal{L} \tilde{u}_r =(\Delta+k^2)\tilde{u}_r \in H^m(\Omega_h)$, there is $f=0$ a.e. in $\Omega$ and 
		\[
		\|f\|_{H^m(\Omega_h)}\leq (2+k^2)\|\tilde{u}_r\|_{H^{m+2}(\Omega_h)}\leq C(\Omega,m,k)\|u_r\|_{H^{m+2}(\Omega)}.
		\]
        Choose $0<h_0<1$ with $2h_0<h$. For $0<\delta<h_0$, take a cut off function $\chi_\delta\in C_c^\infty (\Omega_{2\delta})$ with $\chi_\delta|_{\Omega_\delta}=1$ and $\|D^\alpha\chi_\delta\|_{L^\infty}\leq C_\alpha \delta^{-|\alpha|}$. 

        Since $D^\beta f =0$ a.e. in $\Omega$, applying Lemma \ref{source-extension} with order $m-|\beta|$ yeilds
        \[
        \|D^\beta f\|_{L^2(\Omega_{2\delta})}\leq C(\Omega,m)\delta^{m-|\beta|}\|f\|_{H^{m}(\Omega_h)},\quad |\beta|\leq m.
        \]
        Define $f_\delta=\chi_\delta f$ and extend it by zero outside $\Omega_h$, we have that for $\alpha\leq m$,
        \[
        \|D^\alpha f_\delta\|_{L^2(\mathbb{R}^2)}\leq C(\Omega,m)\sum_{|\beta|\leq |\alpha|}\delta^{-|\alpha|+|\beta|}\delta^{m-|\beta|}\|f\|_{H^m(\Omega_h)}\leq C(\Omega,m)\delta^{m-|\alpha|}\|f\|_{H^m(\Omega_h)}. 
        \]
        Consequently, for $m\geq 1$,
        \[
        \| f_\delta\|_{L^2(\mathbb{R}^2)}\leq C(\Omega,m)\delta^m\|f\|_{H^m(\Omega_h)},\quad \|f_\delta\|_{H^m(\mathbb{R}^2)}\leq C(\Omega,m)\|f\|_{H^m(\Omega_h)}
        \]
        and 
        \(
        \| f_\delta\|_{L^2(\mathbb{R}^2)}= o(1) 
        \)  as $\delta\rightarrow 0$ for $m=0$.
        
        Take a ball $B\Supset \Omega_h $ such that $k^2$ is not a Dirichlet eigenvalue of $-\Delta$ in $B$. Let $w_\delta\in H_0^1(B)$ solve
        \[
        (\Delta+k^2)w_\delta=f_\delta \quad \text{in} \; B.
        \]
        Dirichlet solvability and elliptic regularity in $B$ yield
        \[
        \|w_\delta\|_{H^1(B)}\leq C(B,k)\|f_\delta\|_{L^2(B)},\quad \|w_\delta\|_{H^{m+2}(B)}\leq C(B,m,k)\|f_{\delta}\|_{H^m(B)}.
        \]
        Here $B$ is chosen independent of $\delta$, so is $C(B,m,k)$.

        Take a bounded smooth domain $\Omega_{\delta/2}\Subset\hat{\Omega}\Subset \Omega_\delta$ and set $\hat{u}_r^\delta=\tilde{u}|_{\hat{\Omega}}-w_\delta$, it satisfies 
        %
		%Take $\hat{\Omega}\Subset\Omega^\delta$ with smooth boundary and let $\hat{u}_r^\delta$ be the solution to:
		\[\Delta \hat{u}^\delta_r+k^2 \hat{u}^\delta_r=0,\quad \text{in}\; \hat{\Omega}\]%,\quad \hat{u}_r^\delta |_{\partial \hat{\Omega}}=\tilde{u}_r|_{\partial \hat{\Omega} }.\]
        and
		%
		%
		%\item Approximation estimate
		%
		%
		%Since $Eu_r$ satisfies:
		%\[L(Eu_r)= \hat{f}, \quad \text{on}\; \hat{\Omega}, \quad E{u}_r |_{\partial \hat{\Omega}}=E{u}_r|_{\partial \hat{\Omega} }.\]
		%
		%the difference between $\hat{u}_r^\delta$ and $u_r$ reads:
		%
		%
		\[ 
		\|\hat{u}^\delta_r|_{\Omega}-u_r\|_{H^1(\Omega)} \leq  \|w_\delta\|_{H^1(\hat{\Omega})}  \leq  C(B,k)\|f_\delta\|_{L^2(B)}\leq C(\Omega,k,m)\delta^m\|u_r\|_{H^{m+2}(\Omega)}, \]
		%due to the a priori estimate for $\hat{u}^\delta_r-Eu_r$.
        %
        for $m\geq 1$ while $\|\hat{u}^\delta_r|_{\Omega}-u_r\|_{H^1(\Omega)}=o(1)$ as $\delta\rightarrow 0$ for $m=0$.
        Meanwhile,
        \[
        \|\hat{u}_r^\delta\|_{H^{m+2}(\hat{\Omega})}\leq \|\tilde{u}_r\|_{H^{m+2}(\Omega_h)}+\|w_\delta\|_{H^{m+2}(B)}\leq C(\Omega,m,k)\|u_r\|_{H^{m+2}(\Omega)},
        \]
        for $m\geq 0$, which completes the proof.

\end{proof}

\subsection{Approximation of the Runge part}

By the quantitative approximation proved in \cite{Chen2026}, the Runge part can be approximated by cylindrical-wave modes with spectral accuracy. Moreover, the norm increase of the approximant on the larger domain can be estimated.

\begin{lemma}[\cite{Chen2026}]\label{runge-approx-new}
	Assume that $\Omega$ is a star-shaped domain with smooth boundary. %${\color{blue}\widetilde\Omega}=\{x\in\mathbb{R}^2: dist(x,\Omega)<\delta\}$ ($\delta>0$) and 
	$\Omega\Subset \widetilde\Omega\Subset O_R $ such that $dist(\overline{\Omega},\partial{\widetilde\Omega})>\delta$. $u$ satisfies Helmholtz equation $\Delta u+k^2u=0$ in $\widetilde\Omega$ and there exists $\mu\in L^2(\partial\widetilde\Omega)$ with $\|\mu\|_{L^2(\partial \widetilde\Omega}<M$ such that $u=\int_{\partial\widetilde\Omega}\Phi_k(x,y)\mu(y)\mathrm{d}y$, where $\Phi_k(x,y)$ is the fundamental solution. Assume that $k^2$ is not a Dirichlet eigen value of $-\Delta$ in these domains. Then for any $\varepsilon>0$, $u|_{\Omega}$ can be approximated by a series converging in $O_R$ as
	\[
	u_\varepsilon (x) =\sum_{n\in\mathbb{Z}}c_n J_n(k|x|)e^{\mathrm{i}n\theta},
	\]
	such that
	\[
	\|u|_{\partial\Omega}-u_\varepsilon|_{\partial\Omega}\|_{L^2(\partial\Omega)}\leq \varepsilon M, \quad \| u_\varepsilon  \|_{H^{1/2}(\partial O_R)}\leq C\varepsilon^{-\frac{c}{\delta} }M, 
	\]
	where $c$ and $C$ are constants independent of $u$ and $\delta$.
\end{lemma}

\begin{theorem}[\cite{Chen2026}]\label{runge-part}
	With the assumptions in Lemma \ref{runge-approx-new}, further assume that $ \Omega\Subset \widetilde{\Omega}\Subset O_\rho\Subset O_R$ with $dist(\overline{\Omega},\partial \widetilde\Omega)>\delta$. Then, for sufficiently large $N$, there exist  $\{c_n\}_{n=-N}^N$ such that 
	\begin{equation}\label{num-series}
		u_R^N(x)=\sum_{n=-N}^N c_n \frac{J_n (kr)}{J_n (kR)}e^{\mathrm{i}n\theta}
	\end{equation}
	approximates $u$ with  \[
	\| u-u_R^N\|_{L^2(\partial\Omega)}\leq  \left(\frac{\rho}{R}\right)^{N\cdot \frac{1}{1+\beta}} M,\; \text{and}\quad \left(\sum_{n=-N}^N\sqrt{1+n^2} |c_n|^2\right)^{\frac{1}{2}}\leq C(\Omega,k)  \left(\frac{R}{\rho}\right)^{N\cdot \frac{\beta}{1+\beta}}M,
	\]
	where $\beta=\frac{C}{\delta}$ for some constant independent of $\delta$.
	
	%\[
	%\| h|_{\Omega}-u_N|_{\Omega}  \|_{L^2(\Omega)}\leqslant \varepsilon \|h\|_{H^1(\tilde{\Omega})}+\frac{C}{R^N}\|u_N\|_{L^2(\partial O_R)},\quad \text{and}\quad \|u_N\|_{H^{1/2}(\partial O_R)}\leq C \varepsilon^{-\beta}\|h|_{\Omega}\|_{L^2(\Omega)}.
	%\]
\end{theorem}

\section{Proof of the main results}
\paragraph{Proof of Theorem \ref{main-approximation}}

%{\color{blue} List for the draft, for referring see section 2.}

\newtheorem*{nonumthm}{Theorem}

\iffalse
\begin{theorem}%\label{main-approximation}
	%\color{teal}
	Let \(\Omega\Subset O_R\) be a strictly star-shaped polygon in $\mathbb{R}^2$. $u$ is a solution to the boundary value problem \eqref{eq:Helm_boundary} with $u|_{\partial\Omega}=g\in X_{m,p_0},\, p_0>2, {\color{blue}m\in\mathbb N,\ m>1}
	$.
	
	Then there exist coefficients $c_n$ and {\color{blue}$c_{j,l}$}, such that for $R>1$ with $k^2$ not being Dirichlet eigen value for $O_R$,
	\[
	u_N(x)=\sum_{n=-N}^{N} \frac{c_n}{C_{R,n}}J_n(kr)e^{\mathrm{i}n\theta}  + \sum_{(j,l)\in \mathcal{I}_m} c_{j,l}\sigma_{j,l},\quad C_{R,n}=\frac{1}{|n|!}\left(\frac{kR}{2}\right)^{|n|}	\]
	converges to $u$ as $N\to \infty$, and with an error for $m>1$ as
	\[
	\|u-u_N\|_{L^2(\partial\Omega)} \leq  C \left(\frac{1}{\tilde{R}^{\sqrt{N}}}+\frac{1}{ N^{\frac{m}{2}}}\right)\|g\|_{X_{m,p_0}}.
	\]
	and
	\[
	\|u-u_N\|_{L^\infty(\partial\Omega)} \leq  C \left(e^{-\kappa_1\sqrt{N}}+\frac{1}{ N^{\kappa_2m}}\right)\|g\|_{X_{m,p_0}}.
	\]
	for some $\tilde{R}\in (1,R)$ and $\kappa_1,\kappa_2\in (0,1)$.

	Moreover, $c_n $ have bounds
	\[
	\sum_{n=-N}^N\sqrt{1+n^2} |c_n|^2 \leq C q^{-2N+c\sqrt{N}}\|g\|^2_{X_{m,p_0}},\quad |c_{j,l}|\leq C\|g\|_{X_{m,p_0}}, \quad    \lambda_{l,j} = \frac{l\pi}{\omega_j} \leq m+1. 
	\]
	for some $q\in (0,1)$ and $c>0$, and 
	\[
	q^{2N}\sum_{n=-N}^N\sqrt{1+n^2} |c_n|^2 \leq C e^{-a\sqrt{N}}\|g\|^2_{X_{m,p_0}}.
	\]
\end{theorem}
\fi

\begin{proof}
	By Theorem \ref{singularity-decomposition}, $u(x)=u_r(x)+u_s(x)$, where $u_s$ is the singular part and $u_r\in H^{m+2}(\Omega)$. According to Lemma \ref{harmonic-approx2} , there exists $\hat{u}^\delta_r $ satisfying Helmholtz equation in a smooth domain $\hat{\Omega}$ with \(
	\Omega_{\delta/2}\Subset\hat\Omega\Subset\Omega_\delta\) for $0<\delta<h_0$ such that $\hat u_r^\delta\in H^{m+2}(\hat\Omega)$ and 
	\begin{equation}
		\|\hat{u}^\delta_r|_{\Omega}-u_r\|_{H^1(\Omega)} \leq  \displaystyle C(\Omega,m) \delta^{m} \|g\|_{X_{m,p_0}}, m\geq 1.
	\end{equation}
	Since Lemma \ref{harmonic-approx2} implies %$\hat{u}^\delta_r|_{\Omega}-u_r\in H^{m+2}(\Omega)$,  
	$\|\hat{u}_r^\delta-u_r\|_{H^{m+2}(\Omega)}
	\leq C(\Omega,m,k)\|g\|_{X_{m,p_0}}$,
	by Sobolev interpolation theorem we have that for $s\in [1,m+2]$
	\begin{equation}
		\|\hat{u}^\delta_r|_{\Omega}-u_r\|_{H^s(\Omega)} \leq  \displaystyle C(\Omega,m) \delta^{m\frac{m+2-s}{m+1}  } \|g\|_{X_{m,p_0}}, \quad  m\geq 1.
	\end{equation}
	By Sobolev embedding theorem and trace theorem, 
	\[
	\|\hat{u}^\delta_r-u_r\|_{L^\infty (\partial\Omega)} \leq {\color{blue}C}\|\hat{u}^\delta_r|_{\Omega}-u_r\|_{H^2(\Omega)}   \leq  \displaystyle C(\Omega,m) \delta^{m\frac{m}{m+1}  } \|g\|_{X_{m,p_0}}, m\geq 1.
	\]
	
	Since $\Omega$ is fixed strictly star-shaped polygon, we can choose the smooth outer approximation domains $D_\delta$ and $\hat{\Omega}$ with uniform Lipschitz character and uniform strict star-shapedness constants such that \[\Omega\Subset D_\delta\Subset \hat{\Omega},\quad  dist(\overline{\Omega},\partial D_\delta)\geq c_0\delta,\quad dist(\overline{D}_\delta,\partial\hat{\Omega})\geq c_1\delta,\] 
	where $c_0,c_1$ are constants independent of $\delta$,   by radial smoothing and outward shifting. %{\color{blue}If necessary, choose these domains inside the extension domain supplied by Lemma~\ref{harmonic-approx2} and restrict the extension to the new outer domain; its approximation and Sobolev bounds are preserved.
    By the standard interior and exterior Green representations \cite{McLean2000}, the exterior Dirichlet-to-Neumann estimate \cite{Baskin2016}, and the trace theorem, for the regular part $\hat{u}_r^\delta (x)\in H^{m+2}(\hat{\Omega})$ there exists $\mu\in L^2(\partial\hat{\Omega})$ such that
	\[
	\hat{u}_r^\delta(x)=\int_{\partial\hat{\Omega}}\Phi_k(x,y)\mu(y)\mathrm{d}s_y,\quad x\in\hat{\Omega},
	\]
	with
%	
%	For the regular part $\hat{u}_r^\delta (x)\in H^{m+2}(\hat{\Omega})$ which is bounded by $\|\hat{u}_r^\delta\|_{H^{m+2}(\hat{\Omega})} \leq C(\Omega,m,k)\|g\|_{X_{m,p_0}}$ due to Lemma \ref{harmonic-approx2} and the trace theorem, there exists a single layer potential density $\mu\in H^{m+1/2}(\partial\hat\Omega)$, and 
	\[
	\|\mu\|_{L^2(\partial\hat\Omega)}\leq  C\|\hat{u}^{\delta}_r\|_{H^2(\hat \Omega)} \leq  C(m,p_0,k,\Omega)\|g\|_{X_{m,p_0}}.
	\] 
	The constants are independent of $\delta$, owing to the uniform geometric bounds. 
	By the Rellich estimate and duality, the uniform Lipschitz and strict star-shapedness bounds, together with the uniform Dirichlet spectral gap, imply
	\[
	\|v\|_{L^2(D_\delta)}
	\le C\|v\|_{L^2(\partial D_\delta)},
	\qquad (\Delta+k^2)v=0\ \text{in }D_\delta,
	\]
	with \(C\) independent of $\delta$. The uniform spectral gap follows from nonresonance in $\Omega$ and convergence of the Dirichlet eigenvalues under outer approximation.

	Then, by the quantitative Runge approximation Theorem \ref{runge-part},  we have that there exists an approximant 
	\[u_R^N=\sum_{n=-N}^N c_n \frac{J_n (kr)}{C_{R,n}}e^{\mathrm{i}n\theta},\]
	where
	\[
	c_n=a_n\frac{C_{R,n}}{J_n(kR)},\qquad
	\sup_{n\in\mathbb Z}
	\left|\frac{C_{R,n}}{J_n(kR)}\right|<\infty,
	\]
	with \(a_n\) being the coefficients in Theorem \ref{runge-part} and the boundedness comes from the asymptotic behavior of Bessel function for large $n$ and nonresonance condition on \(O_R\).
	
	$u_R^N$ satisfies Helmholtz equation in $O_R$ that can  approximate $\hat{u}^\delta_r(x)$ in $D_\delta$ with
	\[
	\|\hat{u}_r^\delta (x) - u_R^N(x)\|_{L^\infty(\partial\Omega)} \leq \frac{C}{\delta} \|\hat{u}_r^\delta (x) - u_R^N(x)\|_{L^2(\partial D_{\delta})}\leq  \frac{C}{\delta} q^{N\cdot \frac{1}{1+\beta}}  \|g\|_{X_{m,p_0}}, \quad \beta=\frac{C}{\delta}.
	\]
	for some constant $0<q<1$ which depends on the assumption on the convergence radius of the Runge part.

	Meanwhile, the normalized coefficients
	\[
\left(\sum_{n=-N}^N\sqrt{1+n^2} |c_n|^2\right)^{\frac{1}{2}}\leq C  q^{-N\cdot \frac{\beta}{1+\beta}}\|g\|_{X_{m,p_0}}.
	\]

	Decompose $u$ as
	\begin{align*}
		u(x)& = u_r(x)+u_s(x)\\
		& =   u_r(x) - \hat{u}_r^\delta (x) +\hat{u}_r^\delta (x) - u_R^N(x)+ u_R^N(x) + u_s(x).
	\end{align*}
	%
	%Therefore, by Lemma \ref{runge-part} and \ref{harmonic-approx2} there is
	%
	Then, for $u_N:=u_R^N+u_s$ one has
	\[
	\|u-u_N\|_{L^\infty(\partial \Omega)} \leq \frac{C(k,m,\Omega)}{\delta} q^{N\cdot \frac{1}{1+\beta}}  \|g\|_{X_{m,p_0}}+ C(k,m,\Omega)\delta^{m\frac{m}{m+1}}\|g\|_{X_{m,p_0}}.
	\]
	Since $\beta= \frac{B}{\delta}$ for some constant $B$, there is
	\[
	\frac{1}{1+\beta}\geq C \delta. 
	\]
	Then
	\[
	\|u-u_N\|_{L^\infty(\partial\Omega)} \leq C(k,m,\Omega) \|g\|_{X_{m,p_0}}\left( q^{ CN\delta}\frac{1}{\delta} + \delta^{m\frac{m}{m+1}} \right).
	\]
	Due to the arbitrariness of $\delta$, we can take $\delta = N^{-s}$ and fix $s=1/2$ to get, %and take $\beta_N=B\sqrt{N}$, we have
%	\[
%	\|u-u_N\|_{L^\infty(\partial\Omega)}\leq C(k,m,\Omega)\|g\|_{X_{m,p_0}}\left(\sqrt{N}{\color{blue}q}^{\frac{N}{1+B\sqrt{N}}}+N^{-\frac{m^2}{2(m+1)}}\right)
%	\] 
%	while,
%	\[
%	A_N:= \left(\sum_{n=-N}^N\sqrt{1+n^2} |c_n|^2\right)^{\frac{1}{2}}\leq C  {\color{blue}q}^{-N\cdot \frac{B\sqrt{N}}{1+B\sqrt{N}}}\|g\|_{X_{m,p_0}}.
%	\]
	
	\[
	q^{ CN\delta}\frac{1}{\delta} + \delta^{\frac{m^2}{m+1}}\leq q^{C' N^{1-s}}+\frac{1}{N^{s m\cdot\frac{m}{m+1}}}=:e^{-\kappa_1\sqrt{N}}+\frac{1}{N^{\kappa_2 m}},
	\]
	for $\kappa_2=\frac{m}{2(m+1)} $.
	
	%{\color{blue}For the coefficient bound, $\beta=B\sqrt N$ gives
	%\[
	%-\frac{2N\beta}{1+\beta}=-2N+\frac{2N}{1+B\sqrt N}\ge -2N+c\sqrt N
	%\]
	%for any fixed $0<c<2/B$ and sufficiently large $N$. Squaring the coefficient estimate and using $0<q<1$ yields the stated bound. Multiplication by $q^{2N}$ gives the weighted estimate with $a=c\log(1/q)>0$.}
    
	The $L^2$ case can be treated similarly without using Sobolev interpolation, taking $\delta=N^{-1/2}$,
	\[
	\|u-u_N\|_{L^2(\partial\Omega)}
	\le C\|g\|_{X_{m,p_0}}\left(\delta^m+\delta^{-1}q^{N/(1+\beta)}\right).
	\]
	 
	%If we take 
	%\[
	%\delta=\frac{m}{C\ln\frac{1}{{\color{blue}q}}}\frac{\ln N}{N},\quad  \Rightarrow  {\color{blue}q}^{CN\delta}=\frac{1}{N^{2m}}, \quad \delta^{(m-\frac{3}{4})\frac{m}{m+1}}\leq C\frac{1}{N^{m/2}},\quad m>1.
	%\]

%	\[
%	coef: {\color{blue}q}^{-N\frac{\sqrt{N}C}{1+\sqrt{N}C}},\quad err: {\color{blue}q}^{N\frac{1}{1+\sqrt{N}C}}
%	\]
	
\end{proof}

We remark that this convergence result is not trivial, since the solution is actually approximated by those in larger domains, whereas the solution itself is probably not continuable.

%Remark: there exists $g^{\varepsilon_0}$ 
%\[
%\|g-g^{\varepsilon_0}\|_{L^2(\Gamma)}\leq \varepsilon_0\|g\|_{H^{s}(2\pi)},\quad \|g^{\varepsilon_0}\|\leq C(k,s)\varepsilon_0^{\frac{k-s}{s}}\|g\|_{s}
%\]

\paragraph{Error analysis for the collocation problem}

Approximation of a solution to Helmholtz equation with those on larger domains is ill-posed. Therefore, in the numerical computation of such ill-posed problems, it is important to bridge the gap between discrete and continuous settings.

\begin{theorem}\label{main-approximation2}
	Under the assumptions of Theorem \ref{main-approximation} with $m>1$, let $\{x_i\}_{i=1}^M$ be collocation points on $\partial\Omega$ such that on each side $\Gamma_j$ of $\partial\Omega$, the $M_j$ points form Gauss-Legendre nodes with $M_j\ge cM$. Denote $\mbf{g}=[g(x_1),\cdots g(x_M)]$, and assume that there exists a constant $\mathcal{M}>0$ such that $\|g\|_{X_{m,p_0}}\leq \mathcal{M}$. Denote $\mbf{d}=\{\{c_n\}_{n=-N}^N;\{c_{j,l}\}_{(j,l)\in\mathcal{I}_m} \}$, and $u_N$ as \eqref{solution-form}. 
	Then there exist constants $\kappa_1,\kappa_2,\kappa_3\in (0,1)$ and $\rho_1,\rho_2\in (0,1)$, such that for sufficiently large $M,N$, the minimizer $\mbf{d}^*$ to
	\[
	J(\mbf{d})=\frac{1}{M}\sum_{i=1}^M  |u_N(x_i)-g(x_i)|^2+\alpha \sum_{n=-N}^N |c_n|^2,
	\]
	gives $u_N^*$ with
	\[
	\|u_N^*-g\|_{L^2(\partial\Omega)}\leq C\left(N^{-\kappa_3m}+M^{-\kappa_3m}+ \rho_1^{N^{\kappa_1}}+\rho_2^{ M^{\kappa_2}}\right),
	\]
	provided that $\alpha\sim \left(\tilde{N}^{-\frac{m^2}{2(m+1)}}+\tilde{N}^{\frac{1}{2}}\left(\rho_1^N+\rho_2^M\right)^{\frac{1}{1+\sqrt{\tilde{N}}}}\right)^2\left(\rho_1^{N}+\rho_2^M\right)^{\frac{2\sqrt{\tilde{N}}}{1+\sqrt{\tilde{N}}}}$ where $\tilde{N}=\min\{M,N\}$.
The constant $C$ depends on \(\mathcal{M}, m,p_0,k,\Omega,R\) and the nodes ratio allocation, but not on \(M,N\) and $g$.	
	%$\rho_1>\frac{R_1}{R}$ is that $\Omega\Subset O_{R_1}\Subset {O_R}$, $\rho_2> e^{-4\frac{dist(\Omega,\partial O_R)}{diam(\Omega)} }$.
	\iffalse
	Define the cost functional
	\[
	J(\mbf{d})=\frac{1}{M}\sum_{i=1}^M (u_N(x_i)-g(x_i))^2+\alpha \sum_{n=1}^N C_n^2
	\]
	with $\alpha \sim \left(\frac{\delta^{(m-\frac{1}{2}-\eta)}}{(m-2)!}+(R^{-N}+h)^{\frac{\delta}{1+\delta}}\right)^2\left(R^{-N}+h\right)^{-\frac{2}{1+\delta}} $, where 
	\[
	u_N=\sum_{n=-N}^{N} \frac{C_n}{J_n(kR)}J_n(kr)e^{\mathrm{i}n\theta}  + \sum_{\stackrel{j=1,\cdots,J}{l=1,\cdots, N_j}} c_{j,l}\sigma_{j,l} .
	\]
	Then the minimizer $\mbf{d}^*$ gives the solution $u^*_N$ such that for $\forall \delta>0$,
	\[
	\|u-u_N^*\|_{L^2(\Omega)} \leq  C(k,\mathcal{M},\Omega) \left(\frac{\delta^{(m-\frac{1}{2}-\eta)}}{(m-2)!}+\left(\frac{1}{R^N}+h^M\right)^{C'\delta} \right).
	\]

	where $\{x_i\}$ uniformly distributed points %where $c_i$ is the integral coefficients based on Lagrange interpolation, and  
	and 
	$\alpha \sim \left(\frac{1}{(m-2)!}\left(\frac{1}{N^{\alpha_1}}+\frac{1}{M^{\alpha_1}}\right)+\frac{1}{R^{N^{\alpha_2}} }+h^{M^{\alpha_3}}\right)^2\left(R^{-N}+h^M\right)^{\frac{-2N^\alpha}{1+N^\alpha}} $, 
	then there exist constants $\alpha_1,\alpha_2,\alpha_3\in (0,1)$ such that
	\[
	\|u_N^*-g\|_{L^2(\partial\Omega)}\leq C(k,\mathcal{M},\Omega)\left(\frac{1}{(m-2)!}\left(\frac{1}{N^{\alpha_1}}+\frac{1}{M^{\alpha_1}}\right)+\frac{1}{R^{N^{\alpha_2}} }+h^{M^{\alpha_3}}\right),
	\]
	\fi
\end{theorem}

\,

\begin{proof}
	Since 
	\begin{align*}
		u(x)&=u_r(x)+u_s(x)\\
		&=u_r(x)-\hat{u}^\delta_r(x)+\hat{u}^\delta_r(x)-u_R(x)+u_R(x)-u_R^N(x)+u_R^N(x)+u_s(x).\end{align*}
	Similar to the proof of Theorem \ref{main-approximation}, given $\delta>0, \varepsilon>0$, there exist $c_n'$ and constants $c_{j,k}'$, such that 
	\[
	\tilde{u}_N(x)=\sum_{n=-N}^{N} \frac{c_n'}{C_{R,n}}J_n(kr)e^{\mathrm{i}n\theta}  + \sum_{(j,l)\in \mathcal{I}_m} c_{j,l}'\sigma_{j,l}=:\tilde{u}^N_R+\tilde{u}_s, 
	\]
	approximates $u$ with an error of 
	\begin{equation}\label{err-g-tildeu}
	\|g-\tilde{u}_N\|_{L^\infty(\partial\Omega)}\leq C(\Omega,m)\delta^{\frac{m^2}{m+1}}\|g\|_{X_{m,p_0}}+C(k,\Omega)\frac{\varepsilon}{\delta}\|g\|_{X_{m,p_0}}+C\rho_1^N\|u_R\|_{L^2(\partial O_R)},
	\end{equation}
	for some $\rho_1\in (0,1)$ with 
	\[
	\|u_R\|_{H^{1/2}(\partial O_R)}\leq \varepsilon^{-\beta}\|g\|_{X_{m,p_0}}.
	\]
	where $\beta\sim \frac{C}{\delta}$. As a result,
	\begin{equation}
		\frac{1}{M}\sum_{i=1}^M|\tilde{u}_N(x_i)-g(x_i)|^2\leq C(k,\Omega) \left(\delta^{\frac{m^2}{m+1}}+\frac{\varepsilon}{\delta}+\rho_1^N\varepsilon^{-\beta}\right)^2\|g\|^2_{X_{m,p_0}}=:C(k,\Omega)\varepsilon_0^2\|g\|^2_{X_{m,p_0}} . \label{un0-residual}
	\end{equation}
	
	By definition of $\mbf{d}^*$, $J(\mbf{d}^*)\leq J(\mbf{d}')$. If $\alpha\sim \varepsilon_0^2\cdot \varepsilon^{2\beta} $, then
	\begin{align}
		\frac{1}{M}\sum_{i=1}^M|u_N^*(x_i)-g(x_i)|^2\leq &\frac{1}{M}\sum_{i=1}^M|\tilde{u}_N(x_i)-g(x_i)|^2+\alpha \varepsilon^{-2\beta}\|g\|^2_{X_{m,p_0}}\notag\\
		\leq & C(k,\Omega) \left(\delta^{\frac{m^2}{m+1}}+\frac{\varepsilon}{\delta}+\rho_1^N\varepsilon^{-\beta}\right)^2\|g\|^2_{X_{m,p_0}}\notag\\
		=: & C(\mathcal{M},\Omega)\varepsilon_0^2,\label{un1-residual}
	\end{align}
	Meanwhile,
	\[
	\sum_{n=-N}^N|c_n|^2 \leq \frac{1}{\alpha}C\varepsilon_0^2+ C\varepsilon^{-2\beta}\leq C(\mathcal{M},\Omega) \varepsilon^{-2\beta},
	\]
	which implies that the regular part $\|u_{R}^{*N}\|_{L^2{ (\partial O_R)}}\leq C(\mathcal{M},\Omega)\varepsilon^{-\beta}$, where we have written  $u_N^*=u_{R}^{*N}+u^*_{s}$.

	Let
		\[
		S_m=\operatorname{span}\{\sigma_{j,l}|_{\partial\Omega}:(j,l)\in\mathcal{I}_m\}.
		\]
		Each $\sigma_{j,l}$ admits an edgewise holomorphic extension to fixed side parameter complex neighborhoods $U_i$ of $\Gamma_i$. Since $S_m$ is finite dimensional, Gauss-Legendre quadrature and equivalence of norms imply that for sufficiently large $M$ with $M_i\geq cM$,
		\[
		\max_i\|s\|_{L^\infty(U_i)}\leq C\|s\|_{L^2(\partial \Omega)}\leq C\left(\frac{1}{M}\sum_{\nu=1}^M|s(x_\nu)|^2\right)^{\frac{1}{2}},\quad s\in S_m.
		\]
		Therefore,
		\begin{align*}
		\max_i\|u_{s}^*-\tilde u_s\|_{L^\infty(U_i)}
		&\leq 
		C\|u_{s}^*-\tilde u_s\|_{\ell_M^2}
		\le
		C\left(\|u_N^*-\tilde u_N\|_{\ell_M^2}
		+\|u_{R}^{*N}-\tilde u_R^N\|_{\ell_M^2}\right) \\
		& \le C(\varepsilon_0+\varepsilon^{-\beta})
		\|g\|_{X_{m,p_0}}
		\le C\varepsilon^{-\beta}\|g\|_{X_{m,p_0}},
		\end{align*}
		where \(\|v\|_{\ell_M^2}^2=M^{-1}\sum_{\nu=1}^M|v(x_\nu)|^2\). 
		The regular difference $u_R^{*N}-\tilde{u}_R^N$ also satisfies the same form of bound in these neighborhoods, and so does the total $u_N^*-\tilde{u}_N$.
		
		%For the corner singularity part, for $m\notin\mathbb{Z}$ $\sigma_{j,l}=0$ on $\Gamma_j$ and $\Gamma_{j+1}$ and is analytic on other sides. For $m\in\mathbb{Z}$, the trace of $\sigma_{j,l}$ on the sides are also analytic and analytically continuable on complex plane.   
		%Then $(u_N^*-\tilde{u}_N)|_{\Gamma_i}$ is analytic and can be locally complexified to a holomorphic function that is bounded by $\|u_R^N\|_{L^2(\partial O_R)}+\|u_N^*\|_{L^2(\partial O_R)}$. 
		%
		By Lemma \ref{interpolation} in Appendix A, there exists $\tilde{\rho}_2,\rho_2\in (0,1)$,
		\[
		\|u_N^*-\tilde{u}_N\|_{L^2(\Gamma_i)}
		\le C(\Gamma)\left(\frac1{M_i}\sum_{\nu=1}^{M_i}|u_N^*(x_\nu)-\tilde{u}_N(x_\nu)|^2\right)^{1/2} +
		C(\Gamma,\mathcal{M})\varepsilon^{-\beta}\tilde{\rho}_2^{M_i},
		\]
		and by $M_j\geq cM$
		\[
		\|u_N^*-\tilde{u}_N\|_{L^2(\partial\Omega)}
		\le C(\Omega)\left(\frac1{M}\sum_{\nu=1}^{M}|u_N^*(x_\nu)-\tilde{u}_N(x_\nu)|^2\right)^{1/2} +
		C(\Omega,\mathcal{M})\varepsilon^{-\beta}\rho_2^{M}.
		\]
		Therefore,
		\begin{align*}
			\|u_N^*-g\|_{L^2(\partial\Omega)}&\leq \|u_N^*-\tilde{u}_N\|_{L^2(\partial\Omega)}+\|\tilde{u}_N-g\|_{L^2(\partial\Omega)}\\
			& \leq C(\Omega)\left(\frac1{M}\sum_{j=1}^{M}|u_N^*(x_j)-\tilde{u}_N(x_j)|^2\right)^{\frac{1}{2}} +
			C(\Omega,\mathcal{M})\varepsilon^{-\beta}\rho_2^{M}+\|\tilde{u}_N-g\|_{L^2(\partial\Omega)}%\\
			%& \lesssim  \|\mbf{u}_N^*-\mbf{g}\|_{\mathbb{R}^M}+\|\mbf{g}-\mbf{u}_N\|_{\mathbb{R}^M}+\delta(M)+\|u_N-g\|_{L^2(\partial\Omega)}\\
			%& \lesssim \varepsilon_0+\varepsilon_0+\left(\frac{1}{R}\right)^{c\sqrt{M}}\varepsilon_0^{-\beta}+E(N,m)
		\end{align*}
		By \eqref{un0-residual} and \eqref{un1-residual}
		\begin{align*}
			\frac1{M}\sum_{j=1}^{M}|u_N^*(x_j)-\tilde{u}_N(x_j)|^2 & \leq \frac{2}{M}\sum_{j=1}^{M}|u_N^*(x_j)-g(x_j)|^2+\frac{2}{M}\sum_{j=1}^{M}|g(x_j)-\tilde{u}_N(x_j)|^2\\
			&\leq C(\mathcal{M},\Omega)\left(\varepsilon^2_0+\varepsilon^2_0\right)
		\end{align*}
		%where Lemma \ref{interpolation} in Appendix A is used. 
		Consequently, combining \eqref{err-g-tildeu} yields,
		\[
		\|u_N^*-g\|_{L^2(\partial \Omega)}\leq C(\mathcal{M},\Omega)\left(\delta^{\frac{m^2}{m+1}}+\frac{\varepsilon}{\delta}+\left(\rho_1^N+\rho_2^M\right)\varepsilon^{-\beta}  \right).
		\]
		By arbitrariness of $\delta$ and $\varepsilon$, 
		Take $\delta = C_0\tilde{N}^{-s}$ with $s\in (0,1)$ and $\tilde{N}=\min\{M,N\}$, $C_0$ is sufficiently large that the Runge exponent $\beta\leq \tilde{N}^s$. Enlarging this exponent if necessary, we henceforth use $\beta=\tilde{N}^s$ in the preceding estimates and take 
		\[
		\varepsilon=%\left(\rho_1^{N}+\rho_2^M\right)^{\frac{1}{1+\beta}}\leq 
		\left(\rho_1^{N}+\rho_2^M\right)^{\frac{1}{1+\tilde{N}^s}},
		\] one has 
		\begin{align*}
			&\delta^{\frac{m^2}{m+1}}+\frac{\varepsilon}{\delta}+\left(\rho_1^N+\rho_2^M\right)\varepsilon^{-\beta}\\
			\lesssim & \tilde{N}^{-s\frac{m^2}{m+1}}+C\tilde{N}^s\left(\rho_1^N+\rho_2^M\right)^{\frac{1}{1+\tilde{N}^s}}+\left(\rho_1^N+\rho_2^M\right)^{\frac{1}{1+\tilde{N}^s}}\\
			\lesssim &\tilde{N}^{-s\frac{m^2}{m+1}}+\tilde{N}^{s}\left(\rho_1^N+\rho_2^M\right)^{C\tilde{N}^{-s}}\\
			\lesssim & N^{-\kappa_3m}+M^{-\kappa_3m}+ \rho_1^{N^{\kappa_1}}+\rho_2^{ M^{\kappa_2}}
		\end{align*}
		with $\kappa_1,\kappa_2,\kappa_3\in (0,1)$ being some constants. %and $m$ depending on the number of singular samples. 
		 %
		%Combining the estimate for $\|\tilde{u}_N-g\|_{L^2(\partial\Omega)}$ in Theorem \ref{main-approximation}, and 
		Taking $s=\frac{1}{2}$ in the error estimate and in $\alpha\asymp\varepsilon_0^2\varepsilon^{2\beta}$ gives the asserted conclusion and parameter choice.

\end{proof}

\section{Numerical examples}
The first numerical example considers L-shape geometry. % with which we want to illustrate the effect of additional designed data solution. 
There is geometry singularity at the convex corner. We consider the solution to 
\begin{equation*}
\left\{
\begin{array}{ll}
\triangle u + k^2u = -5-\frac{5k^2}{4}(x^2+y^2)\,, & \text{in}\quad  \Omega\\
u=0 & \text{on} \quad \partial\Omega.
\end{array}
\right.
\end{equation*}
Due to the existence of the corner, there are components of the source the not in the image of $\{u\in W_2^{p}(\Omega), u|_{\partial\Omega}=0\}$ through the operator, i.e., the shift theorem for a smooth domain no longer applicable which limit the global regularity. By subtracting the solution $u_0(x,y) = -\frac{5}{4}(x^2+y^2)$, it turns to solving the homogeneous equation with inhomogeneous boundary. We take $k=0.3$ in this example.

\begin{figure}[H]
    \centering
    \includegraphics[width=0.4\linewidth]{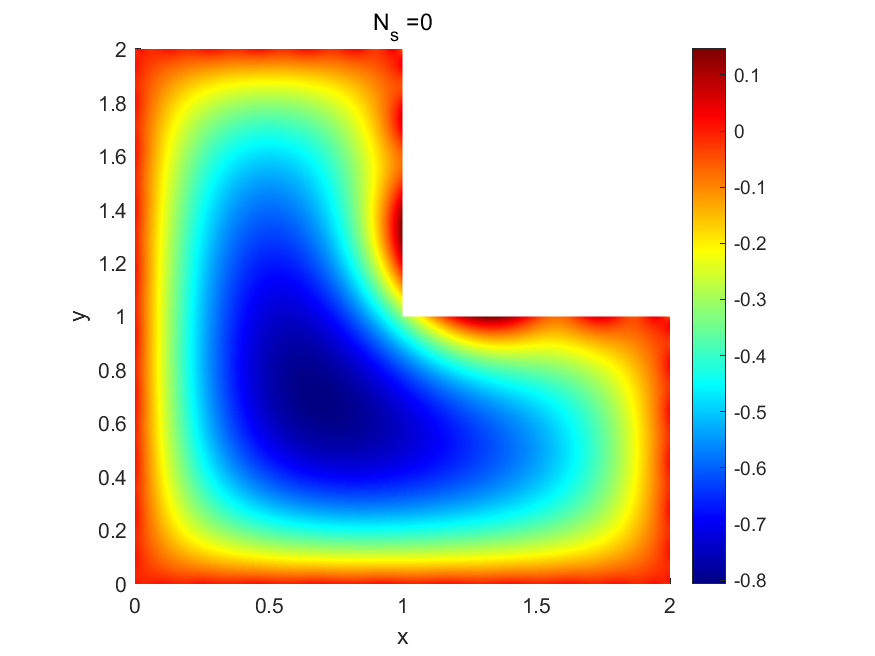}
    \includegraphics[width=0.4\linewidth]{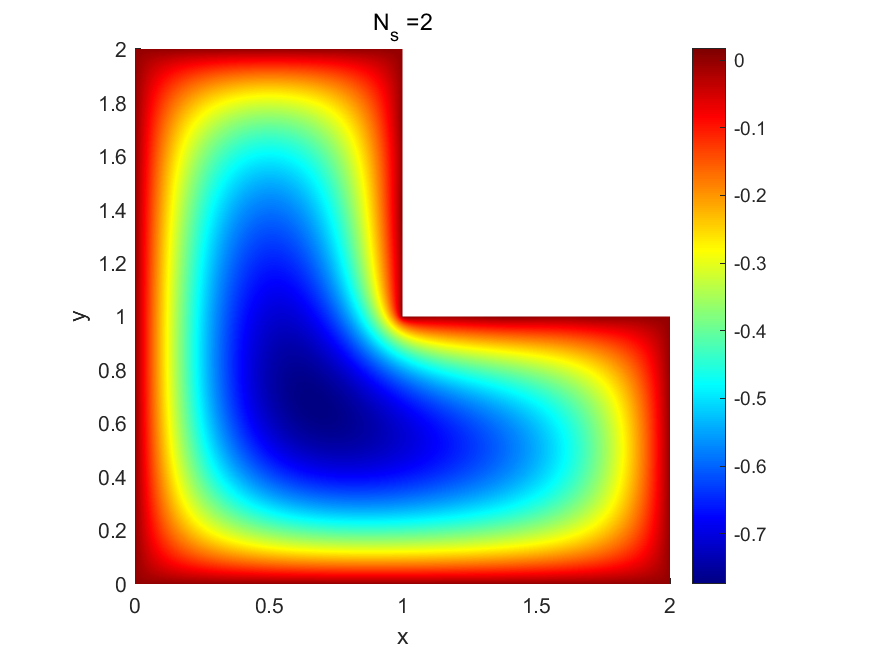}
    \includegraphics[width=0.4\linewidth]{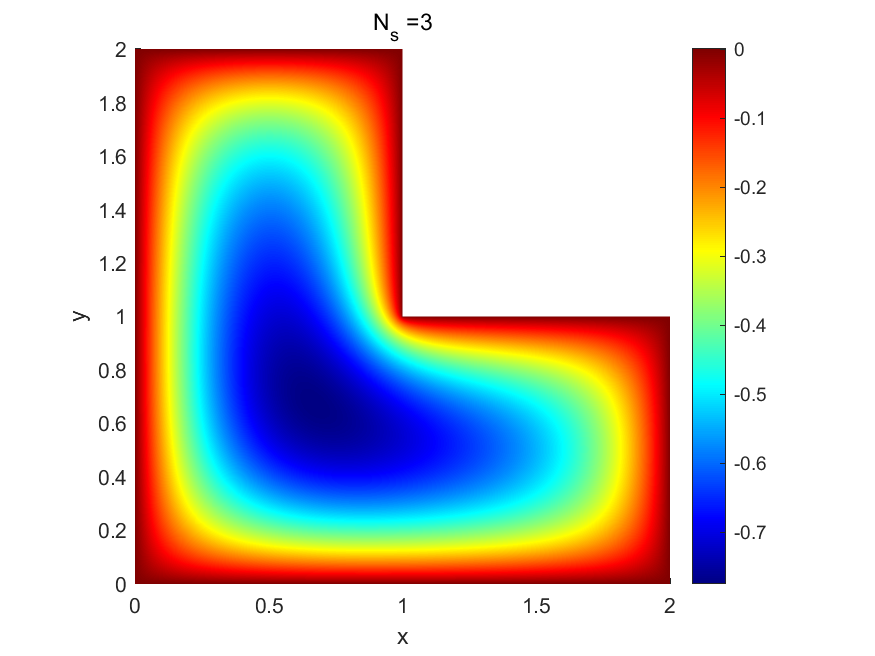}
     \includegraphics[width=0.4\linewidth]{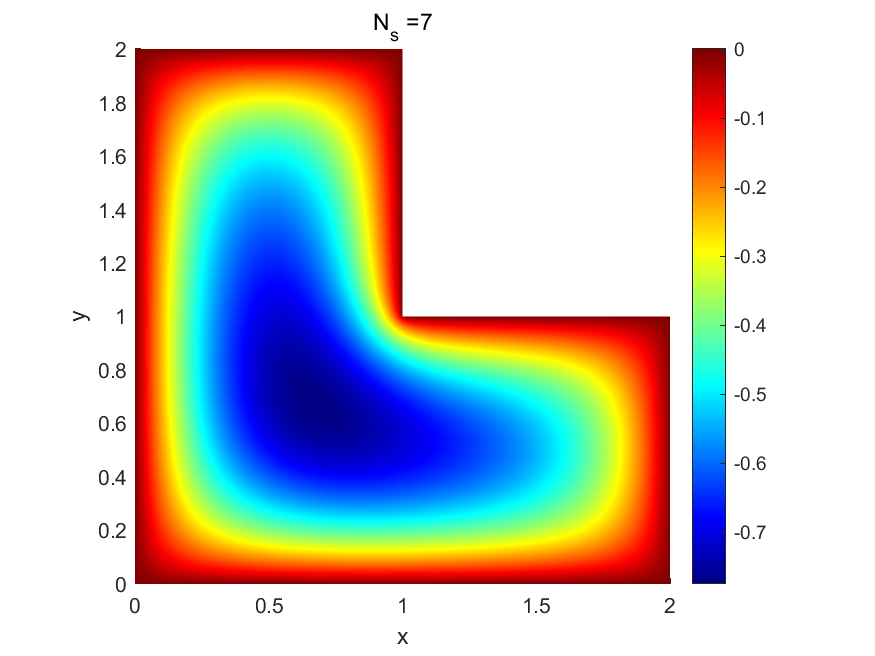}
    \caption{Numerical results for the pentagon example, with various number of singularity samples are included in the approximation.}
        \label{fig:ex1-1}
\end{figure}

\begin{figure}[H]
    \centering
    \includegraphics[width=0.4\linewidth]{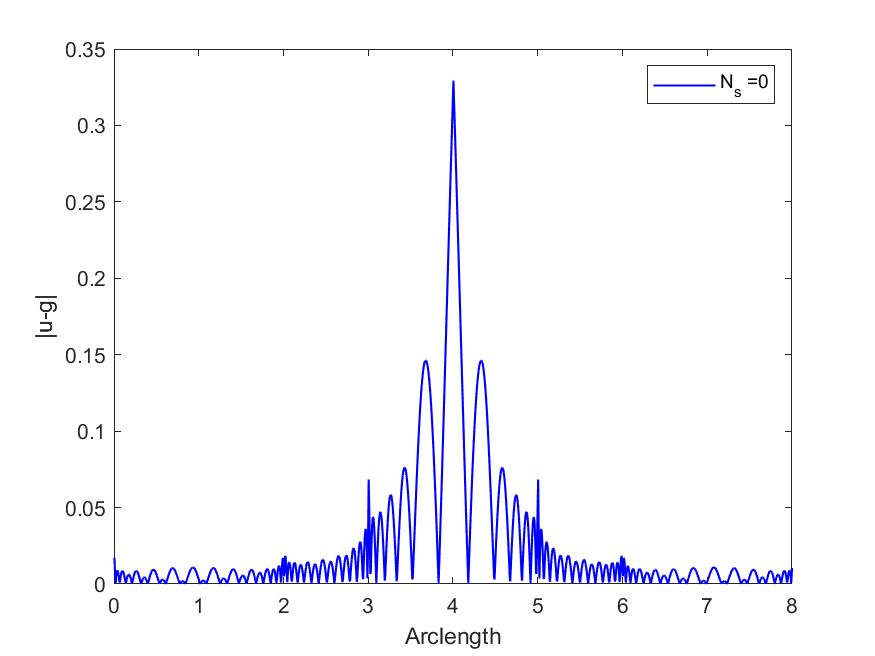}
    \includegraphics[width=0.4\linewidth]{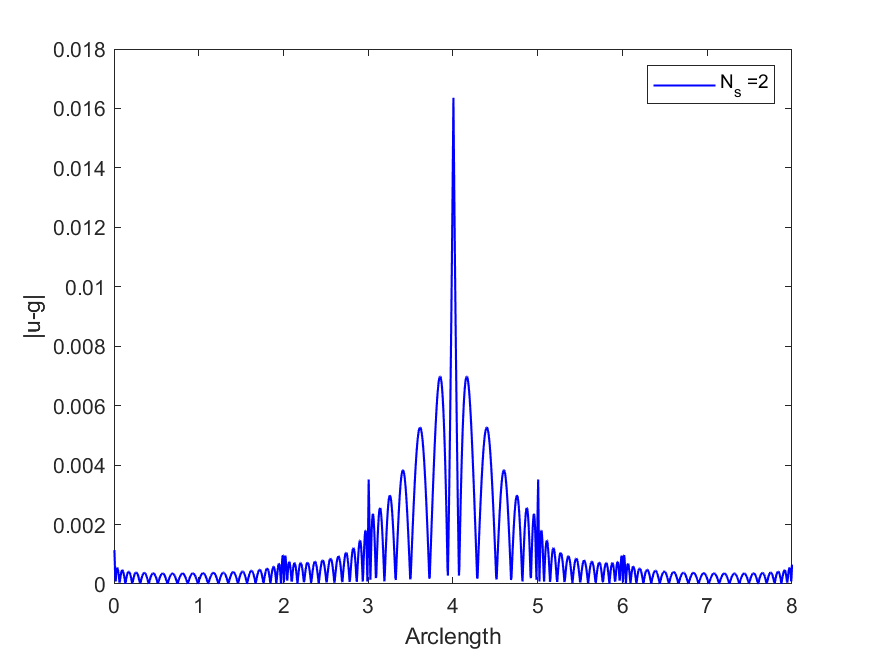}
    \includegraphics[width=0.4\linewidth]{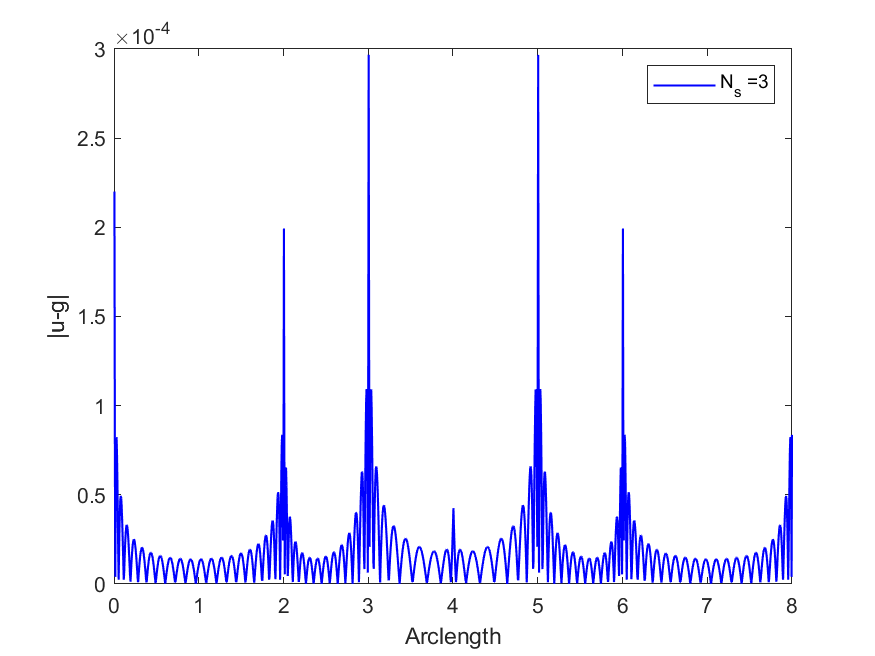}
     \includegraphics[width=0.4\linewidth]{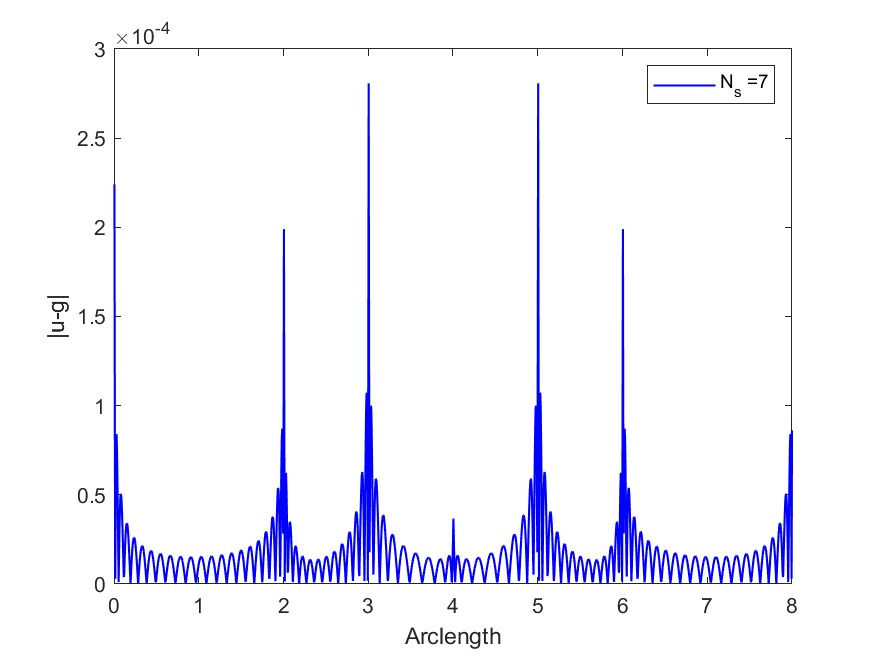}
    \caption{The absolute error distribution along the boundary against arclength, corresponding to the cases in Fig. \ref{fig:ex1-1} }
        \label{fig:ex1-2}
\end{figure}

Fig \ref{fig:ex1-1} gives the solutions by the approximation method in Theorem \ref{main-approximation2} with $0$, $2$, $3$ and $7$ corder singular samples $\sigma_{j,l}$ included. We take $N=40$ and 1600 boundary points are used for the collocation. It can be seen that without corner samples, the result is significantly different from other cases. To further illustrate the accuracy, we plot the boundary residual distributions against the boundary arc-length in Fig. \ref{fig:ex1-2} (the corresponding non-homogeneous boundary value problem ($g=u_0|_{\partial\Omega}$). The case without singular sample indeed shows a large boundary error around the concave corner ($s=4$), even with so many collocation points. Incorporating two singular samples reduces the error from $0.35$ to about $0.016$. Similar results for the case with one singular sample, which lifts the remaining part to $H^2(\Omega)$ and agrees with the previous theoretical analysis. The third sample corresponds to the integer eigen value $\lambda_{j,l}=2$, and the logarithmic type of singularity near the corner. Once this sample is added, the error near this corner is reduced so that it is no longer dominant. And there is no significant change in the error distribution with more singular samples. 

The error near convex corners can be reduced by introducing the corresponding singular samples similarly. Since $\lambda_{j,l}=l\pi/\omega_j=4l$ ($\omega_j=\pi/4$) are integers, only higher order samples are involved and the regularity/error level without such samples is already acceptable ($10^{-4}$). 

We remark that applying [\cite{Kuttler1978}, Eq. (7)], the boundary $L^2$ error controls the interior error of the solution as
\[
\|u-u_\varepsilon\|_{L^2(\Omega)}\leqslant \frac{C_\Omega}{d}\|u-u_\varepsilon\|_{L^2(\partial\Omega)}
\]
where $d:=\min_j |k^2-E_j|/E_j$, $E_j$ are the Dirichlet eigenvalues of the domain and $C_\Omega$ is a domain-dependent constant.
%}

%plot the boundary error by (without, order 2, order 3, order 7) to illustrate the necessity of the log type solution.

%Second example gives the the pentagon, (order 2,7, boundary plot)

\begin{figure}[H]
    \centering
    \includegraphics[width=0.4\linewidth]{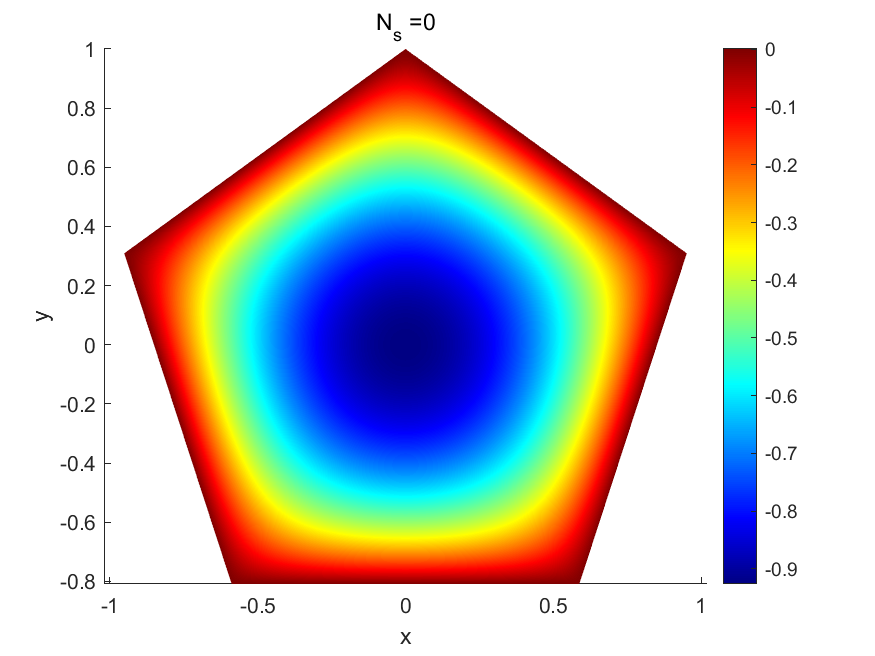}
    \includegraphics[width=0.4\linewidth]{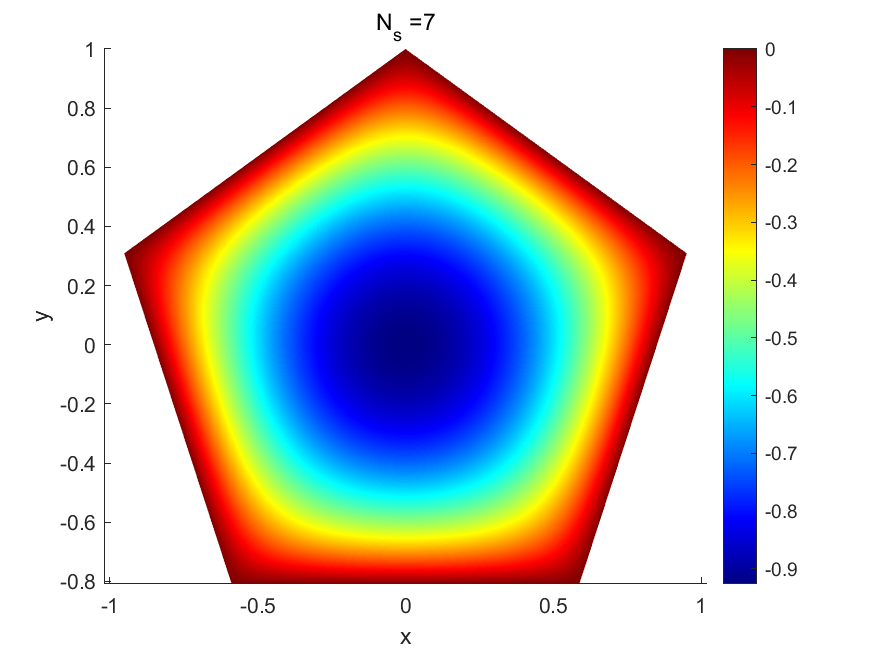}
    \includegraphics[width=0.4\linewidth]{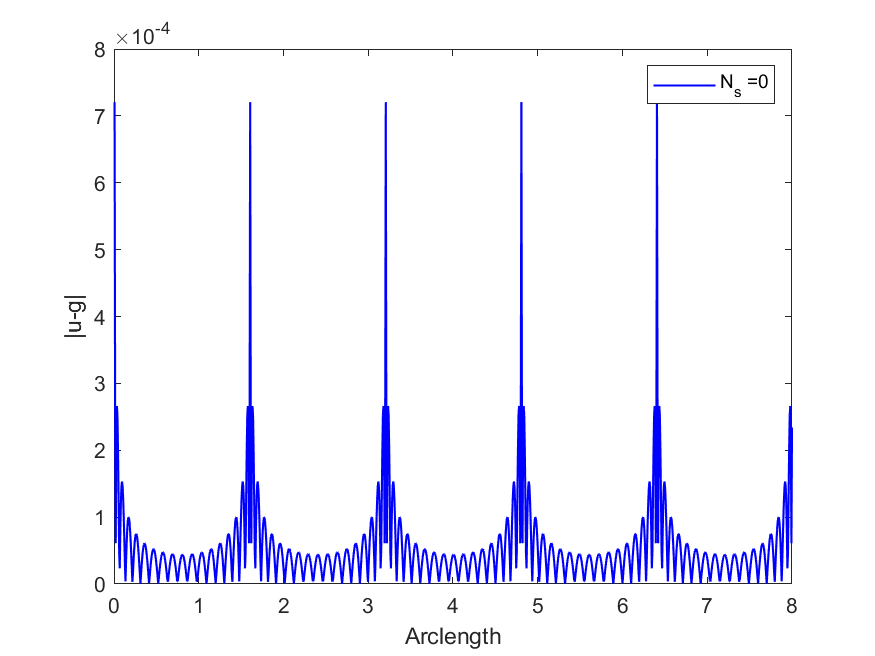}
     \includegraphics[width=0.4\linewidth]{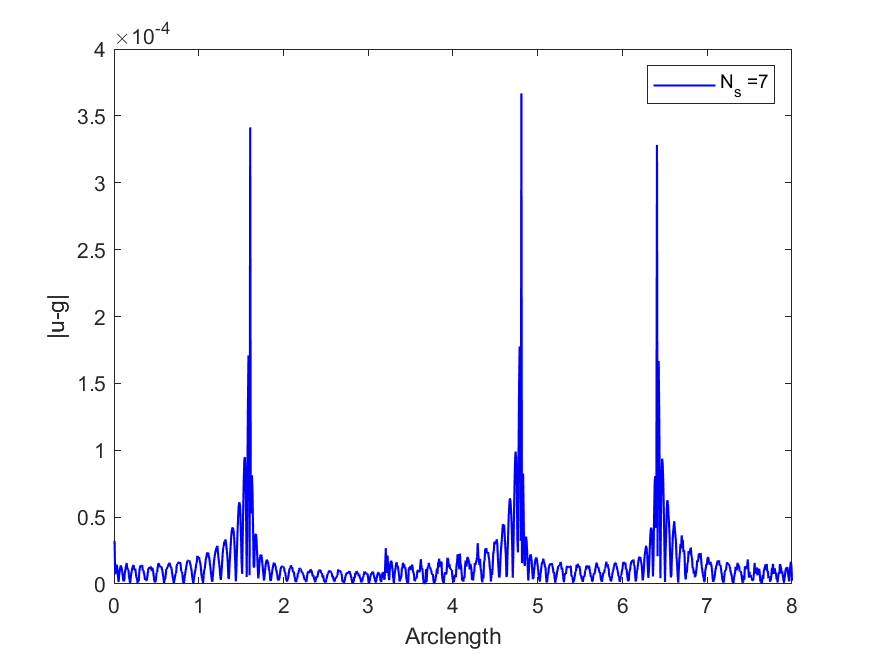}
    \caption{Numerical results for the pentagon example, with and $k=0.3$}
        \label{fig:ex2-1}
\end{figure}

\begin{figure}[H]
    \centering
    \includegraphics[width=0.4\linewidth]{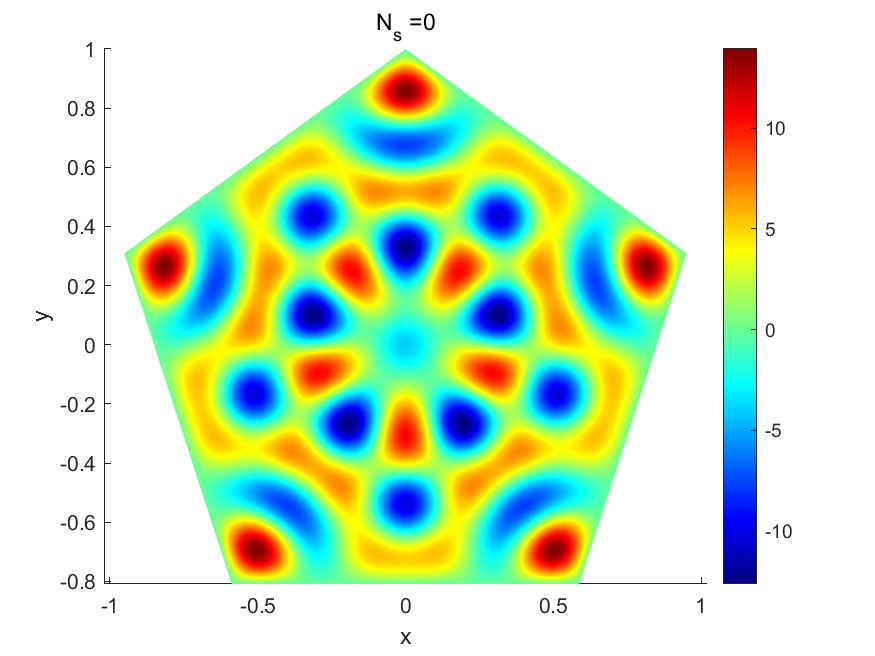}
    \includegraphics[width=0.4\linewidth]{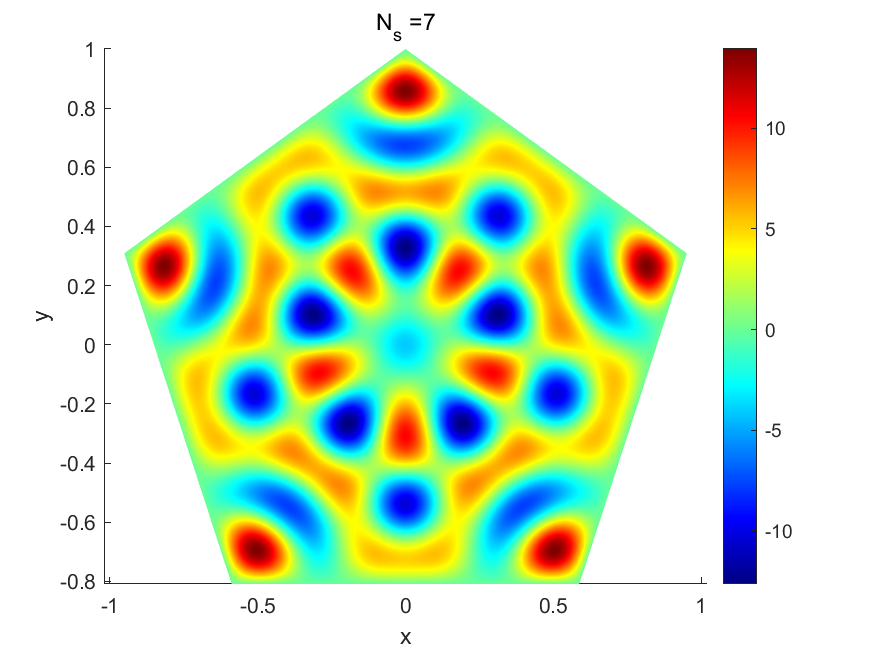}
    \includegraphics[width=0.4\linewidth]{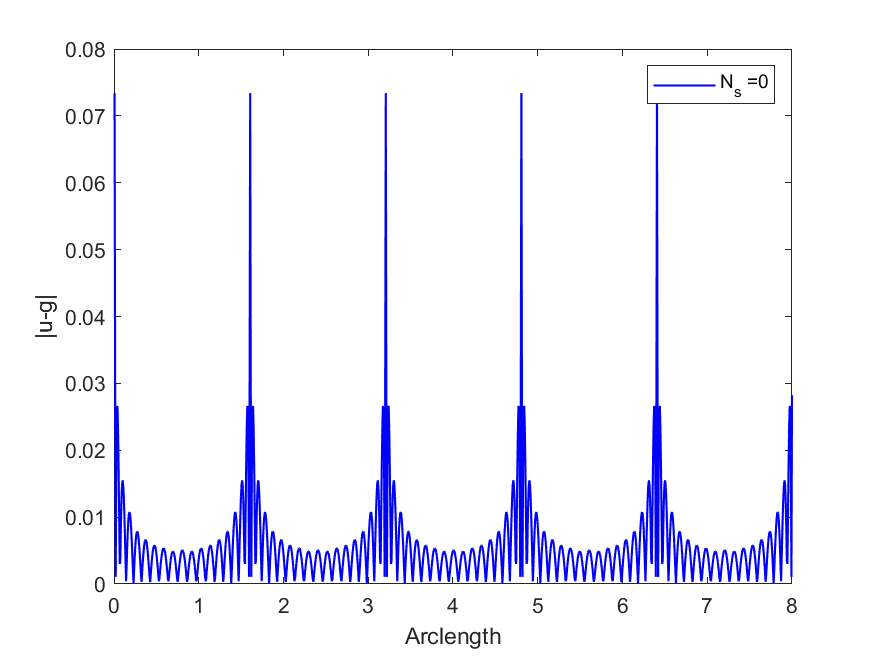}
     \includegraphics[width=0.4\linewidth]{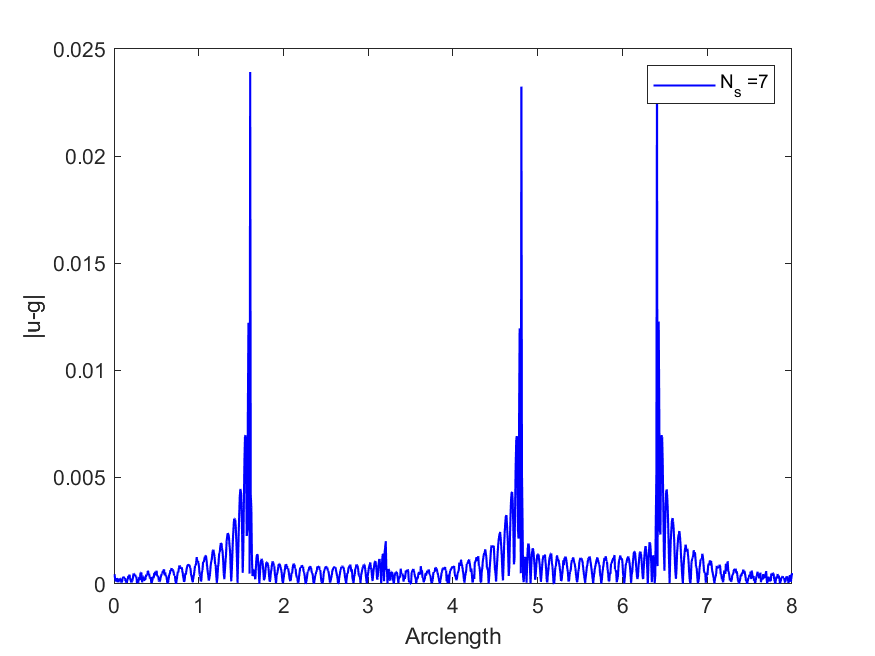}
    \caption{Numerical results for the pentagon example, with and $k=20$}
    \label{fig:ex2-2}
\end{figure}

The second example considers the pentagon shape. The source is the same as the previous case. Fig. \ref{fig:ex2-1} and \ref{fig:ex2-2} give the results with $k=0.3$ and $k=20$ respectively. Since $\omega_j=3/5\pi$ and the solution is naturally in $H^2(\Omega)$ with the present boundary value, we can see that the corner singularity issue is not as severe as the first example. Even so, the subtraction of singular components can further reduce the residual on the boundary (see the difference between $N_s=0$ and $N_s=7$). 

As the wavenumber increases, the norm of the solution increases, so does the absolute boundary residual. The effect of inclusion of singularity samples are similar.

\section{Concluding remarks}
%A summary of the highlights of the work ...

In this work, we give a quantitative approximation for solutions to Helmholtz equations in non-smooth domains by cylindrical wave modes and corner singular solutions. The regularity analysis enables singular decomposition of the solution space. Based on this, the regular part can be approximated by the Runge part, which is further quantified by the explicit Runge approximation. Combining a priori choice of regularization parameters and discrete error, we finally get the convergence rate for continuous and discrete problems. Due to the equivalence of the approximation problem and learning based method, the results provide error estimate which ensures interpretability and generalizability of the choice of learning samples and the reconstruction method.

We remark that the choice of singular samples is compatible with the work of Trefethren \cite{Gopal2019}. In fact, the singular solution
\[
J_{\lambda}(kr)e^{\mathrm{i\lambda\theta}}
\]
can be extended circumferentially and satisfies the Helmholtz equation in the whole plane except for a ray. Thus, it can be represented by single layer potential along the ray, i.e., approximated by poles along the ray. A detailed analysis is worth further study.
\begin{figure}[H]
\centering
\includegraphics[width=0.7\textwidth]{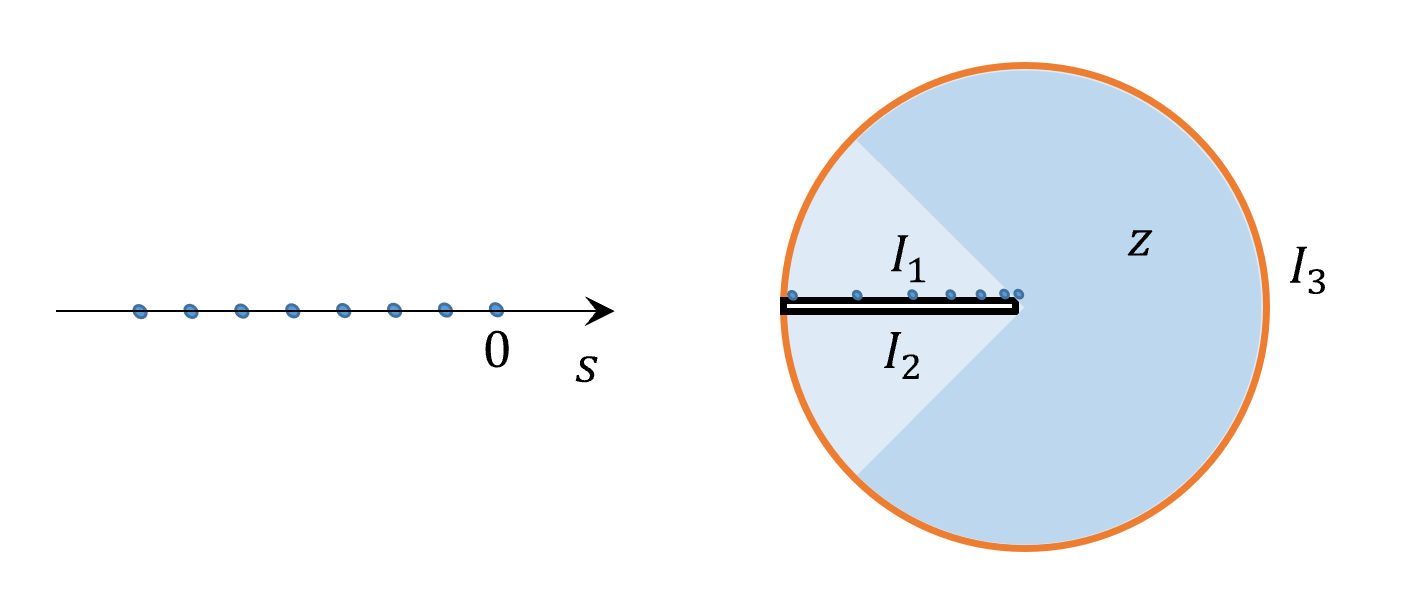}
\caption{The illustration of branch-point type solution and approximation by clustered poles.}
\end{figure}

Besides, the regular/Runge part, i.e., the cylindrical wave modes can be replaced by other samples in larger domains, such as fundamental solutions. 

Similar ideas can be adopted in the analysis of other types of boundary conditions. The present results can be extended to the case of an exterior problem that is often raised in scattering problems. This will be discussed in our forthcoming paper.

\section*{Acknowledgments}
This work was supported by China's National Key Research and Development Programs (No. 2024YFA1012401), China's National Natural Science Foundation (Nos. 12241103, 12571458), Program for Innovative Research Team of Shanghai University of Finance and Economics (No. 2022110917), and funds from Laboratory of Mathematics for Nonlinear Science, Fudan University.

The authors thank Prof. Min Zhong from Southeast University and Mr. Lifu Song from Fudan University for helpful discussions.

\addcontentsline{toc}{section}{

\newpage

\appendix

\section{Appendix I}

%\newpage
\begin{lemma}[discrete control on an analytic segment)]\label{interpolation}
	
	Let \(\gamma\) be a line segment of length \(L\). Suppose that \(f|_\gamma\) admits a holomorphic extension to a fixed complex neighborhood of \(\gamma\) and satisfies
	\[
	|f(z)|\le B
	\]there. Let \(x_1,\ldots,x_N\) be the \(N\)-point Gauss–Legendre nodes on \(\gamma\). Then there exist constants \(C>0\) and \(\rho\in (0,1)\), independent of \(f\) and \(N\), such that
	\[
	\|f\|_{L^2(\gamma)}
	\le C\sqrt L\left(\frac1N\sum_{j=1}^{N}|f(x_j)|^2\right)^{1/2} +
	CB\sqrt L\,\rho^{N}.
	\]
	The constants \(C\) and \(\rho\) depend only on the size of the complex neighborhood of \(\gamma\).

\end{lemma}

\begin{proof}
	Let \(w_j\) be the normalized Gauss–Legendre weights. Parameterize the sides by a real interval. Applying analytic Gauss quadrature \cite{Trefethen2019} (Theorem 19.3) to
	\[
	F(z)=f(z)\overline{f(\overline z)}
	\]gives
	\[
	\left|
	\frac1L\int_\gamma |f|^2\,ds
	-\sum_{j=1}^{N}w_j|f(x_j)|^2
	\right|
	\le CB^2\rho^{2N}.
	\]The classical estimate for Gauss–Legendre weights \cite{Gautschi2004},
	\[
	0<w_j\le \frac{C}{N},
	\]then yields
	\[
	\sum_{j=1}^{N}w_j|f(x_j)|^2
	\le
	\frac{C}{N}\sum_{j=1}^{N}|f(x_j)|^2,
	\]which proves the result.
\end{proof}

%\addcontentsline{toc}{chapter}{参考文献}
%\bibliographystyle{elsarticle-num}%{unsrt}%{plain}
%\bibliographystyle{unsrt}
%\bibliographystyle{elsarticle-num}
%\bibliography{mybibfile-chen-pub-2}
\end{document}